\documentclass{amsart}
\usepackage
{hyperref}

\usepackage{amsmath,amsthm,amssymb,verbatim,mathrsfs,tikz-cd}%
\usepackage{color}
\usepackage{url}

\newtheorem{thm}{Theorem}[section]
\newtheorem{cor}[thm]{Corollary}
\newtheorem{lem}[thm]{Lemma}
\newtheorem{prop}[thm]{Proposition}
\newtheorem{con}[thm]{Conjecture}
\newtheorem{hyp}[thm]{Hypothesis}

\theoremstyle{definition}

\theoremstyle{definition}

\newtheorem{rem}[thm]{Remark}

\newcommand\be{\begin{equation}}
\newcommand\ee{\end{equation}}
\newcommand\bee{\begin{equation*}}
\newcommand\eee{\end{equation*}}
\newcommand\ben{\begin{enumerate}}
\newcommand\een{\end{enumerate}}

\def\G{\ensuremath {{\mathbb G}}}
\def\CC{\ensuremath {{\mathbb C}}}
\def\H{\ensuremath {{\mathcal H}}}
\def\I{\ensuremath {{\mathcal I}}}

\def\R{\ensuremath {{\mathbb R}}}
\def\F{\ensuremath {{\mathcal F}}}

\def\A{\ensuremath {{\mathbb A}}}

\def\P{\ensuremath {{\mathscr P}}}
\def\E{\ensuremath {\mathcal E}}
\def\L{\ensuremath {\mathscr L}}
\def\Z{\ensuremath {{\mathbb Z}}}
\def\D{\ensuremath {\mathcal D}}
\def\C{\ensuremath {\mathscr C}}

\def\K{\ensuremath {\mathcal K}}
\def\G{\ensuremath {\mathcal G}}

\def\T{\ensuremath {\mathscr T}}
\def\SI{\ensuremath {S\mathcal{I}}}

\def\g{\ensuremath {\mathfrak g}}

\def\a{\ensuremath {\mathfrak a}}

\def\GL{\ensuremath {\mathrm{GL}}}

\def\bs{\ensuremath {\backslash}}

\def\temp{\ensuremath {\mathrm{temp}}}
\def\el{\ensuremath {\mathrm{ell}}}

\def\disc{\ensuremath {\mathrm{disc}}}

\def\reg{\ensuremath {\mathrm{reg}}}

\def\ss{\ensuremath {\mathrm{ss}}}
\def\der{\ensuremath {\mathrm{der}}}

\def\cusp{\ensuremath {\mathrm{cusp}}}

\numberwithin{equation}{section}

\title{On the definition of stable transfer factors}
\author{Tian An Wong}
\email{tiananw@umich.edu}
\address{University of Michigan, Dearborn MI, USA}

\begin{document}
\begin{abstract}
We construct stable geometric and spectral transfer factors for a general reductive group and develop some of their basic properties, assuming the refined local Langlands correspondence. Using our definition of stable geometric transfer factors, we show that the stable transfer conjecture for orbital integrals implies the stable transfer of characters and vice versa. The latter is also implied by local Langlands, and in particular this establishes archimedean stable geometric transfer. Finally, we show how the stable geometric transfer factors can be used to define stable spectral transfer factors. 
\end{abstract}

\maketitle

\tableofcontents \addtocontents{toc}{\protect\setcounter{tocdepth}{1}}

\section{Introduction}


\subsection{The stable transfer of distributions}
\label{stintro}

Let $G$ be a connected reductive group over a local field $F$ of characteristic zero. The notion of stable transfer was first formulated by Langlands in \cite{ST}, primarily in the context of SL(2). By stable transfer, we really mean one of two things: On the spectral side, given a local $L$-parameter $\phi$ of $G(F)$ and suitable test function $f$, we write $f^G(\phi)$ for the associated stable character of $f$ at $\phi$. If $\phi'$ is a local $L$-parameter on another group $G'(F)$ that transfers to $\phi$, then the desired transfer $f'$ should satisfy an identity of stable characters
\[
f^{G'}(\phi') = f^G(\phi).
\]
This we call {\em stable spectral transfer}. Here we are identifying the test function $f'$ with the stable character $f^{G'}$, though $f' = f^{G'}$ should be determined only up to its stable character.

On the geometric side, given a strongly regular conjugacy class $\delta$ of $G(F)$, we write $f^G(\delta)$ for the stable orbital integral of $f$ at $\delta$. 
Then we also expect that the spectral transfer induces a transfer of stable orbital integrals,
\[
f^G\to f^{G'} 
\]
hence the overlap of notation, which is standard following Arthur. This we call {\em stable geometric transfer}. Once again  $f' = f^{G'}$ should be determined only up to its stable orbital integral. Implicit in the notation is the expectation that the stable geometric transfer and stable spectral transfers agree with each other.

Whereas in endoscopy, the endoscopic geometric transfer $f^e= f^{G^e}$ was given in terms of the endoscopic transfer factor, also called the Langlands-Shelstad transfer, roughly
\[
f^e(\delta^e) = \sum_{\gamma} \Delta(\delta^e,\gamma) f^G(\gamma),
\]
it is expected that the stable geometric transfer $f'$ requires a stable transfer factor 
\[
f'(\delta') = \int_{\Delta(G/Z)} \Theta(\delta',\delta)f^G(\delta) d\delta,
\]
where now $\Theta(\delta',\delta)$ is a stable distribution on $G'$ and $G$ that is often referred to as the {\em stable transfer factor}, and $\Delta(G)$ is the set of stable regular semisimple conjugacy classes of $G$.

These stable transfer factors have only been constructed in special cases \cite{ST,J,JL}, while Shelstad has sketched a general approach in the archimedean case \cite{shelstad}. In particular, we note that in \cite{ST} and \cite{shelstad}, only in special cases of SL(2) is the stable transfer factor constructed explicitly, and in \cite{J,JL} for special cases related to GL(2). In related work, Sakellaridis has advanced a theory of transfer operators in the relative setting, c.f. \cite{Sak}, we are motivated here by the so-called group case.

\subsection{This paper}

We give a general formulation of stable geometric and spectral transfer factors for general quasisplit connected reductive groups $G$ over a local field $F$. We will  assume the refined local Langlands correspondence, for the most part, where by refined we mean that the correspondence is uniquely characterised by endoscopic character identities \cite{kal}.  We show that the latter implies the existence of these transfer factors abstractly, and we shall also propose an explicit formula for the geometric transfer. We prove properties related to it that, if the case of endoscopy is any indication, will be needed for the primitisation of the trace formula. Our construction builds directly on Arthur's works, so that it is readily adapted for application to the stable trace formula in full generality.

More precisely, we first associate a mesoscopic datum $(G',\G',\xi')$ with auxiliary datum $(\tilde G',\tilde\xi')$, defined in Section \ref{std} and following, which can be viewed as a weakened or `beyond' endoscopic datum. In general, one expects a stable transfer 
\[
f \to f'=f^{\tilde{G}'}
\]
of suitable test functions on $G$ to the space of stable orbital integrals on $\tilde G'$ satisfying the stable character identity 
\begin{equation}
\label{bemap}
f'(\phi') = f^G(\tilde\xi'\circ\phi'),
\end{equation}
where $\phi'$ is a bounded Langlands parameter for $\tilde G'$ and $f'(\phi')$ is the stable character of $f'$ at $\phi'$ (resp. $f^G(\phi)$). This can be rephrased in terms of the existence of a function $f'$ on $\tilde{G}'$ such that the character identity holds, where the function $f'$ is not uniquely determined but its stable orbital integral is. Taking the refined local Langlands conjecture for $G,G'$ as known (due to Shelstad in the archimedean case \cite{sheltrans}), it is possible to specify the transfer in terms of Paley-Wiener functions on the space of tempered Langlands parameters on $G$ and $G'$, denoted $\Phi(G)$ and $\Phi(G')$ respectively. The stable transfer for functions $f$ is then a consequence of the refined local Langlands conjecture, stated as Proposition \ref{thm1}. Note that we impose a simplifying Hypothesis \ref{hyp} that $\xi'(\hat G')$ and $\hat G$ have equal rank, which can be removed with a more careful analysis by stratifying the transfer map according to Levi subgroups. 

As noted by Arthur \cite[\S11]{arthurwork}, and stated differently by Langlands in \cite[\S2]{ST}, the stable spectral transfer should lead to the transfer of stable orbital integrals. That is, the transfer $f'$ is required to satisfy 
\be
\label{bemap2}
f'(\delta') = \int_{\Delta(G/Z)} \Theta(\delta',\delta)f^G(\delta) d\delta,
\ee
where the distribution $\Theta(\delta',\delta)$ is the focus of our study. Here $f^G(\delta)$ is the stable orbital integral of $f$ at the stable conjugacy class $\delta$. The existence of the stable transfer factor follows as an application of the Schwartz kernel theorem and the stable spectral transfer (Corollary \ref{skt}).
For applications to the trace formula, it will likely be necessary to have explicit spectral and geometric characterisations of the transfer $f'$. 

\subsection{Stable geometric transfer factors} 
The key result of this paper is an explicit formula for $\Theta(\delta',\delta)$, and we distinguish our construction with a subscript $\Theta_{\tilde\xi'}(\delta',\delta)$. Our stable geometric transfer factors rely heavily on the local character relations satisfied by stable orbital integrals and stable characters developed by Arthur in \cite{LCR}, which we recall in Section \ref{gtf}. In particular, we show that they satisfy an adjoint relation similar to their unstable counterparts, a result that appears to be new. 
Our geometric transfer factors, defined in Section \ref{sgtf}, take the form of a distribution
\[
\Theta_{\tilde\xi'}(\delta',\delta) =\int_{\Phi(\tilde G',\tilde\zeta')}S'(\delta',\phi')S(\tilde\xi'\circ\phi',\delta)d\phi',
\]
where $S'(\delta',\phi')$ and $S(\phi,\delta)$ are the stable character and the kernel of the Fourier inversion of stable orbital integrals respectively. As an illustration, consider the Fourier transform and its inversion
\[
f(x) = \int_{-\infty}^{\infty}\int_{-\infty}^{\infty} e^{2\pi i(x-y)\xi} f(y) dy \ d\xi,
\]
formally interchanging the order of integration, the inner integral becomes the delta distribution. The distribution $\Theta_{\tilde\xi'}(\delta',\delta)$ is essentially a generalisation of this observation. As we are working with stable objects, much of the theory of endoscopy is essential to our constructions. Our definition of the stable transfer factor is modeled after Shelstad's heurstics \cite{shelstad}. We first formulate in Conjecture \ref{maincon} that $\Theta_{\xi'}$ realizes the kernel $\Theta$ of \eqref{bemap2}, then prove it at the end of Section \ref{stransfer}.

\subsection{Stable spectral transfer factors} 

With our stable geometric transfer factor in hand, in Section \ref{spectf} we apply the preceding results to construct stable spectral transfer factors. Importantly, their definition depends on the stable geometric transfer factors and the surjectivity of a certain stable transfer map $\T^\F$ in \eqref{TF}, which we prove unconditionally in the archimedean case, and conditional on local Langlands in the nonarchimedean case. Assuming these, we show that the spectral transfer factors $\Theta_{\tilde\xi'}(\phi',\phi)$ are well-defined.
 
  It is clear that our results in this paper rely heavily on local Langlands. It is an important question to ask how one might give an intrinsic definition of the stable geometric transfer factor $\Theta_{\tilde\xi'}(\delta',\delta)$ without recourse to it. (Though it should be noted that to even state the characterization of the geometric transfer requires local Langlands.) For the moment, we only explore the surjectivity of $\T^\F$ without the use of local Langlands. We show in Appendix \ref{appb} that one can instead turn to a study of the descent of stable geometric transfer factors along the lines of Langlands and Shelstad in the endoscopic case \cite{LS}, which we initiate but do not complete in Section \ref{descsec}, and also a stable analogue \eqref{wkernel} of Waldspurger's kernel formula for Fourier transforms on Lie algebras \cite{transfert}. Assuming these two identities instead, we can then also define the stable spectral transfer factors.

The broader goal of course, is the primitisation of the stable trace formula, following Arthur's formalisation of Langlands' beyond endoscopy proposal. As we hope to show in future work, this framework of stable transfer will lay the local foundations for beginning work on the problem of primitisation as outlined in \cite{problems}.

\section{Mesoscopic data: A simplified case}
\label{data}

We begin first by fixing notation and recalling basic definitions in \S2.1--2.3, and then introduce the notion of meoscopic data and related constructions in \S2.4--2.6. 

\subsection{Preliminaries}

Let $F$ be a local field of characteristic zero, with an algebraic closure $\bar{F}$. Let $G$ be a quasisplit connected reductive group over $F$. We denote by $\L(M)$ the collection of Levi subgroups of $G$ containing $M$, $\L^0(M)$ the subset of proper Levi subgroups in $\L(M)$, and $\P(M)$ the collection of parabolic subgroups of $G$ containing $M$. Let $A_M$ be the maximal split torus in the centre of a Levi subgroup $M$ of $G$. We then identify the Weyl group of $(G,A_M)$ with the quotient of the normaliser of $M$ by $M$, thus
\[
W(M)=W^G(M)= \text{Norm}_G(M)/M.
\]
If $M_0$ is a minimal Levi subgroup of $G$, which we shall assume to be fixed, and denote $\L=\L(M_0), \P=\P(M_0), \L^0=\L^0(M_0),$ and $W^G_0=W^G(M_0)$. Also write $P_0$ for a minimal parabolic (i.e., Borel) subgroup containing $M_0$. Also, we fix a maximal compact subgroup $K$ of $G(F)$, which is hyperspecial when $F$ is nonarchimedean and $G$ unramified over $F$. 

As usual, we form the real vector space $\a_M= \text{Hom}(X(M)_F,\R)$ where $X(G)_F$ is the module of $F$-rational characters on $G$. We note by $H_G:G(F)\to \a_G$ canonical homomorphism defined by 
$
e^{\langle H_G(x),\chi\rangle} = |\chi(x)|,$ for $x\in G(F),\chi\in X(G)_F,
$
where $|\cdot|$ is the normalised valuation on $F$. We set $\a_{G,F}= H_G(G(F))$ and $\tilde\a_{G,F}=H_G(A_G(F))$, which are closed subgroups of $\a_G$, with associated dual spaces $\a^\vee_{G,F}=\text{Hom}(\a_{G,F},2\pi i \Z)$ and $\tilde\a^\vee_{G,F}=\text{Hom}(\tilde\a_{G,F},2\pi i \Z)$, which are closed subgroups of $i\a^*_G$.  If $F$ is nonarchimedean, all four groups are lattices; if $F$ is archimedean, we have $\tilde{\a}_{G,F} = \a_{G,F} = \a_G$ and $\tilde\a^\vee_{G,F}= \a^\vee_{G,F}=\{0\}$. Fixing a Haar measure on $\a_G$, we obtain a dual Haar measure on the real vector space $i\a^*_G$. If $F$ is nonarchimedean, we normalise measures so that $\a_G/\tilde\a_{G,F}$ and $i\a^*_G/\tilde\a^\vee_{G,F}$ have volume 1. It follows that the volume of the quotient $i\a^*_{G,F} = i\a^*_G/\a^\vee_{G,F}$ equals the index $|\a_{G,F}/\tilde\a_{G,F}|.$

Let $\Gamma=\Gamma_F$ and $W_F$ be the Galois and Weil groups of $\bar{F}/F$ respectively. Let $G^*$ be a quasi-split inner form of $G$ with inner twist $\psi:G\to G^*$. In other words, $\psi$ is an isomorphism such that $\psi\circ\sigma(\psi)^{-1}$ is an inner automorphism for all $\sigma\in\Gamma$. Moreover, fix a bijection of canonical based root data $\Psi(G)^\vee\to \Psi(\hat{G})$, where $\hat{G}$ is the complex dual group of $G$.  Let $(B,T)$ be a Borel pair of $G$ where $B$ is a Borel subgroup of $G$ and $T$ a maximal torus of $B$, not necessarily defined over $F$. For each pair $(B,T)$ in $G$ and $(B_1,T_1)$ in $\hat{G}$ we have a canonical isomorphism $\hat{T}\simeq T_1$.  Define a pinning by $\E=(B,T,\{X_\alpha\})$ where  $\{X_\alpha\}$ runs over simple roots $\alpha$ of $T$ acting on the Lie algebra of the unipotent radical of $B$, and $X_\alpha$ is an element of the eigenspace associated to $\alpha$. Any two pinnings are related by the adjoint action ad$_g$ for some $g\in G_\text{sc}$, the simply connected cover of $G$, unique up to translation by the centre $Z(G_\text{sc})$. The restriction of ad$_g$ to $B$ and $T$ are uniquely determined, so we may define a canonical pinning by taking the inductive limit over all pinnings of $(B,T)$. If we denote by $\rho(w)$ the action of $\Gamma$ on $\hat G$, we can define a new action $\rho_g(w) = \text{ad}_g\rho(g)\text{ad}_{g^{-1}}$ that fixes the original pinning, giving an exact sequence
$
W_F\to \Gamma\to \text{Out}(\hat{G}).
$
Then $^LG$ is isomorphic to $\hat{G}\rtimes W_F$ under this action, sending $(h,w)$ to $(h\rho(w)(g)g^{-1},w)$. The $L$-isomorphism $^L\psi:{^LG} \to {^LG^*}$ induced by $\psi$ allows us to identify the two groups. (We recall that an $L$-homomorphism here is a continuous homomorphism that is analytic on $\hat G$, semsimple on $W_F$ in the sense that the image of any $w\in W_F$ in $\hat G^*$ is semisimple, and commutes with the projections onto $W_F$.)

\subsection{$K$-groups}
We shall work with $K$-groups, following \cite[\S2]{ArtTW}, which streamlines endoscopy theory over archimedean local fields. It is an algebraic variety constructed in the following manner: If $F$ is $p$-adic, then $G$ is just an ordinary connected reductive group, whereas if $F$ is archimedean, then $G$ can have several connected components
\[
G  = \coprod_{\alpha} G_\alpha\qquad \alpha\in \pi_0(G),
\]
a variety whose connected components $G_\alpha$ are reductive groups over $F$, equipped with an equivalence class of frames $(\psi,u) = \{(\psi_{\alpha\beta},u_{\alpha\beta}): \alpha,\beta\in\pi_0(G)\}$ satisfying natural compatibility conditions given in \cite[\S2]{ArtTW}. Here $\psi_{\alpha\beta}: G_\alpha\to G_\beta$ is an isomorphism over $\bar{F}$, and $u_{\alpha\beta}$ is a locally constant function from $\Gamma = \text{Gal}(\bar{F}/F)$ to the simply connected cover $G_{\alpha,\text{sc}}$ of the derived group of $G_\alpha$. Any connected reductive group is a component of a $K$-group that is unique up to weak isomorphism. 

We call $G^*$ a quasisplit inner twist of $G$ if $G^*$ is a connected, quasisplit group over $F$ equipped with a $G^*$-inner class of compatible inner twists $\psi_\alpha:G_\alpha\to G^*$ and a corresponding family of compatible, locally constant functions $u_\alpha:\Gamma\to G^*_{\alpha,\text{sc}}$. We shall call a $K$-group $G$ quasisplit if one of the isomorphisms $\psi_\alpha$ is defined over $F$, i.e., $G_\alpha$ is quasisplit. Unless otherwise indicated, we shall assume $G$ to be a quasisplit $K$-group over $F$.

The usual definitions for connected groups extend to $K$-groups in a natural way. For example, there are similar notions of a Levi $K$-subgroup $M$ of $G$, with associated sets $\L(M)$ and $\P(M)$. The isomorphism $\psi_{\alpha\beta}$ induces a bijection of Borel pinnings from $G_\alpha$ to $G_\beta$, and taking inverse limits their canonical pinnings are thus equivalent and Galois equivariant. A central induced torus $Z$ of a $K$-group $G$ will have central embeddings $Z\simeq Z_\alpha\subset Z(G_\alpha)$ for each $\alpha$, where $Z(G_\alpha)$ is the center of $G_\alpha$, and $\zeta$ determines a character $\zeta_\alpha$ for each $\alpha$. 
These isomorphisms are required to be compatible with the isomorphisms $\psi_{\alpha\beta}$ and $\psi_\alpha$ respectively. We shall call such a pair $(Z,\zeta)$ a central datum for $G$. Finally, we note that $G(F)/Z(F) = G/Z(F)$.

\subsection{Stable conjugacy classes}

Let $G'$ be a reductive group with an embedding of semisimple conjugacy classes into that of $G$. We recall that a semisimple element $\gamma'\in G'(\bar{F})$ is called $G$-{regular} if the image of its conjugacy class in $G(\bar{F})$ consists of regular semisimple elements, and {strongly $G$-regular} if the image consists of strongly regular elements, that is, whose centralisers are tori. If $c$ is a semisimple conjugacy class of $G$, we write $G_{c,+}$ for the centraliser of a representative of $c$ in the component $G_\alpha$ that contains $c$, and write $G_c$ for the identity component of $G_{c,+}$.  We call $c$ elliptic if it lies in an elliptic maximal torus in $G_\alpha$ modulo the split component $A_G\simeq A_{G_\alpha}$ of the centre of $G$.

We write $\Gamma_\ss(G)=\Gamma_\ss(G(F))$ for the set of semisimple conjugacy classes of $G(F)$, $\Gamma(G)= \Gamma_\reg(G(F))$ for the subset of strongly regular, semisimple conjugacy classes in $G(F)$, and $\Gamma_\el(G) = \Gamma_{\reg,\el}(G(F))$ for the subset of regular elliptic conjugacy classes. That is, 
\[
\Gamma_{\el}(G)\subset \Gamma_\reg(G) \subset \Gamma_\ss(G).
\]
We also write $\Gamma_G(G')=\Gamma_{G\text{-reg}}(G'(F)) $ and $\Gamma_{G,\el}(G')=\Gamma_{G\text{-reg},\el}(G'(F))$ for the set of $G$-regular (resp. $G$-regular elliptic) conjugacy classes in $G'(F)$. Clearly each of these sets are equal to the disjoint union over $\alpha$ of the corresponding sets for each connected component, e.g., 
\[
\Gamma(G) = \coprod_{\alpha\in\pi_0(G)}\Gamma(G_\alpha)
\]
and so on. The Weyl group $W(M)\simeq \prod_{\alpha}\text{Norm}_{G_\alpha}(M_\alpha)/M_\alpha$ acts on $\Gamma_{G,\el}(M)$, and we have a decomposition
\[
\Gamma(G) = \bigoplus_{\{M\}}\Gamma_{G,\text{ell}}(M)/W(M),
\]
where the direct sum ranges over conjugacy classes of $K$-Levi subgroups $M$ in the sense of \cite[p.221]{ArtTW}.

We say that two semisimple elements $c_1\in G_{\alpha_1}$ and $c_2\in G_{\alpha_2}$ are stably conjugate if there is a $g_1\in G_{\alpha_1}(\bar{F})$ such that the mapping 
\[
\varphi =\text{Int}(g_1)\circ\psi_{\alpha_1\alpha_2}:G_{\alpha_2}\to G_{\alpha_1}
\]
 maps $c_2$ to $c_1$, and has the property that for any $\sigma\in\Gamma$, the automorphism $\varphi\circ\sigma(\varphi)^{-1}$ of $G_{c_1}$ is inner.  Let $\Delta_\ss(G)=\Delta_\ss(G(F))$ be the set of semisimple stable conjugacy classes in $G(F)$. There is a canonical injective mapping $\delta\to \delta^*$ from $\Delta_\ss(G)$ to $\Delta_\ss(G^*)$, which is a bijection if $G$ is quasisplit. We also define subsets \[
\Delta_{\el}(G)\subset \Delta_\reg(G)=\Delta(G)\subset \Delta_\ss(G)
\]
as above, $\Delta_G(G')=\Delta_{G\text{-reg}}(G'(F)) $ and $\Delta_{G,\el}(G')=\Delta_{G\text{-reg},\el}(G'(F))$ analogously. We also define 
\be
\label{dell}
\Delta(G) = \bigoplus_{\{M\}}\Delta_{G,\el}(M)/W(M)
\ee
for the $W(M)$-orbits as above. For any maximal torus $T$ of $G$ over $F$, we have the finite abelian group $\K(T) = \pi_0(\hat{T}^\Gamma/Z(\hat{G})^\Gamma)$. Given $\gamma\in\Gamma(G)$, there is a bijection between the set of $G(F)$ classes in the stable conjugacy class $\delta$ of $\gamma$ and the set of characters on the group $\K_\delta=\K_\gamma=\K(G_\gamma)$, so we set $n(\delta) = |\K_\delta|.$ (When $F$ is archimedean, this is only true because we are taking $G$ to be a $K$-group.)

\subsection{Mesoscopic datum}
\label{std}

We caution that as was with endoscopy, it it is likely necessary to refine the datum introduced here, which is simply a weakened version of endoscopic datum. Let us call a {\em mesoscopic datum} for a connective reductive group $G$ over $F$, or mesoscopic datum for short, a tuple $(G',\G',\xi')$, where 
\ben
\item 
$G'$ is a connected quasisplit reductive group over $F$,
\item
$\G'$ is a split extension 
\[
1 \to \hat G' \to \G'\to W_F \to 1
\]
such that the homomorphism $W_F\to \text{Out}(\G')$ given by this extension coincides with the homomorphism $W_F\to\text{Out}(\hat G')$,
\item
$\xi'$ is an admissible $L$-embedding of $\G'$ into $^LG$. 
\een
As usual, we shall use $G'$ to stand in for the triple itself. Denote by $G'_{\xi'}$ the connected quasisplit reductive group whose dual group is equal to $\hat{G}_{\xi'} = \xi'(\G') \cap \hat{G}$.  In contrast to the endoscopic setting, it will be important to distinguish between the groups $G'$ and $G'_{\xi'}$. 
We say that a mesoscopic datum $G'$ is \emph{elliptic} if 
the connected component of $\Gamma$-invariants of the centres satisfy $(Z(\hat{G}'_{\xi'})^\Gamma)^0 = (Z(\hat{G})^\Gamma)^0$. The latter condition is also equivalent to the property that 
\[
\left|\text{Cent}(\xi'(\G'),\hat{G})/Z(\hat{G})^\Gamma\right|<\infty,
\]
where $\text{Cent}(\xi'(\G'),\hat{G})$ is the centralizer of $\xi'(\G')$ in $\hat G$, and $Z(\hat G)$ the centre of $\hat G$. 

Given mesoscopic data $(G',\G',\xi')$ and $(G'_1,\G'_1,\xi'_1)$, we say they are isomorphic if there exist an $F$-isomorphism $\alpha:G'_1\to G'$, an $L$-isomorphism $\beta:\G'\to\G'_1$, and an element $g\in\hat{G}'_1$ such that $\alpha:\Psi(G'_1) \to \Psi(G')$ and $\beta:\Psi(\hat{G}') \to \Psi(\hat{G}'_1)$ are dual, and $\text{Int}(g) \circ \xi_1' \circ \beta = \xi'$. We denote by $\F_\el(G)$ the set of isomorphism classes of elliptic mesoscopic datum for $G$, and let $G'$ be a representative of such a class. Let $\text{Aut}_G(G')$ be the set of $g\in\hat{G}$ that induce an isomorphism of mesoscopic data, i.e., Int$(g)$ is an $L$-isomorphism of $\G'$ onto itself, so that we may identify the group of outer automorphisms as 
$
\text{Out}_G(G') \simeq \text{Aut}_G(G')/\xi'(\hat{G}').
$
Any element in $\text{Out}_G(G')$ can be identified with an outer automorphism of $G'$ which is defined over $F$. 

The definitions extend to $K$-groups in a straightforward manner. If $G'$ is a mesoscopic datum for the component $G_\alpha$, then so it is also for $G_\beta$ for any $\beta\in\pi_0(G)$. We can therefore view $G'$ as a mesoscopic datum for the $K$-group $G$. We shall write $\F(G)$ for the set of isomorphism classes of mesoscopic data $G'$ for $G$ that are {\em relevant} to $G$, by which we mean that there is an element in $\Delta_{G\text{-reg}}(G')$ that is an image of some element in $\Delta_\reg(G)$ under $\xi'$, in the sense defined below in Section \ref{ss}. We also write $\F_\el(G)$ for the set of elliptic mesoscopic data.

\subsection{Images of semisimple elements}
\label{ss}

The isomorphism $\hat{T}\simeq T_1$ sends the coroots of $T$ in $G$ to the roots of $T_1$ in $\hat{G}$, the $B$-simple coroots to the $B_1$-simple roots, and the Weyl group of $T$ with contragredient action to the Weyl group of $T_1$. We impose the following simplifying assumption for convenience:  
\begin{hyp}\label{hyp}
Let $G'$ be a mesoscopic datum such that that $\xi'(\hat G')$ and $\hat G$ have equal rank.
\end{hyp}
If $(B_1',T_1')$ is a Borel pair in $\hat G'$ then there is an $x\in\hat{G}$ such that $\text{Int}(x)\circ\xi'$ maps $T'_1$ to $T_1$ and $B'_1$ to $B_1$. If $(B',T')$ is a Borel pair in $G'$, then we have an isomorphism $\xi'(\hat{T}')\simeq \hat{T}$ defined by the composition
\[
\xi'(\hat{T}')\to \xi'(T_1')\to T_1 \to \hat{T}.
\]
Let us write $T'_{\xi'}$ for the torus dual to $\xi'(\hat T')$, so that  $T'_{\xi'} \subset G'_{\xi'}$. Dualizing the maps above gives an isomorphism $T'_{\xi'}\simeq T$. These isomorphisms map the coroots of $T'_{\xi'}$ in $G'_{\xi'}$ to a subsystem of coroots of $T$ in $G$, the Weyl group $W_{T'_{\xi'}}$ of $T'$ into a subgroup of the Weyl group $W_T$ of $T$, and the roots of $T'_{\xi'}$ into a subset of the roots of $T$. Then the map 
$
T'_{\xi'}/W_{T'_{\xi'}}\to T/W_{T}
$
of Weyl group orbits is independent of all choices. Since these orbits classify the conjugacy classes of semisimple elements in $G'_{\xi'}(\bar{F})$ and $G(\bar{F})$, we thus have a canonical map of semisimple conjugacy classes from $G'_{\xi'}(\bar{F})$, and in fact $G'(\bar F)$, to $G(\bar{F})$. 

Suppose that $T'$ is defined over $F$, and recall that we are assuming that $G$ is quasisplit. If $(B,T)$ is Borel pair in $G$ such that $T$ and $T'_{\xi'} \to T$ are defined over $F$ following Steinberg's theorem, then we say $T'_{\xi'}\to T$ is an {\em admissible embedding} of $T'_{\xi'}$ in $G$. It is determined up to conjugation by elements in
\[
\{g \in G_\text{sc}(\bar{F}): g\sigma(g^{-1})\in T(\bar{F}), \sigma\in\Gamma\},
\]
where as usual $G_\text{sc}$ denotes the simply connected cover of the derived group of $G$. We say a strongly $G$-regular $\gamma' \in G'(F)$ is an {\em image of} $\gamma\in G(F)$ if $\gamma$ lies in the image of the stable conjugacy class of $\gamma'$. For arbitrary $G$-regular semisimple $\gamma'$, we set $T' = G'_{\gamma'}$ and choose an admissible embedding of $T'$ in $G$. Then if $\gamma\in G(F)$ is regular semsimple and $T = G_\gamma$ then we say that $\gamma'$ is an {\em image of} $\gamma$ if there exists $x\in G$ such that $\text{Int}(x)$ maps $\gamma$ to the image $\gamma$ of $\gamma'$ under $T'\to T$. The correspondence $(\gamma',\gamma)$  is  independent of the choice of admissible embedding. 

\begin{rem}
If $G'$ has the same rank as $G$, it follows from Borel-de Siebenthal theory regarding maximal reductive subgroups of complex reductive groups that $G'$ is simply an elliptic endoscopic group of $G$. Yet even in this case, the problem of stable transfer remains to be addressed unless $G'$ is isomorphic to $G$. On the other hand, if we do not require the ranks of $G'_{\xi'}$ and $G'$ to be equal, we have to stratify the stable transfer mapping according to Weyl orbits of Levi subgroups of $G$ in a manner similar to \cite{LCR}, but 
we shall avoid this situation for the time being. Note that Shelstad's work suggests that it is also possible to reduce to the equal rank case \cite{shelstad}. 
\end{rem}

\subsection{Auxiliary data}

The group $\G'$ need not be an $L$-group, so there might not be an $L$-isomorphism from $\G'$ to $^LG'$ which is the identity on $\hat{G}'$. Thus given any $G'\in\F(G)$, we shall fix an auxiliary datum $(\tilde{G}',\tilde\xi')$ where $\tilde{G}'\to G'$ is a $z$-extension, by which we mean a split central extension of $G'$ by an induced torus $\tilde{C}'$, and $\tilde\xi':\G'\to {^L\tilde{G}'}$ is an $L$-embedding satisfying the conditions of \cite[Lemma 2.1]{LCR}. Namely, we require that the $z$-extension 
\[
1\to \tilde{C}'\to \tilde{G}'\stackrel{r}{\to} G'\to1
\]
over $F$ satisfies:
\ben
\item
the central subgroup $\tilde{C}'$ is an induced torus,
\item
the dual exact sequence $1\to \hat{G}'\to\hat{\tilde{G}}'\to \hat{\tilde{C}}'\to1$ extends to a short exact sequence of $L$-homomorphisms
\[
1\to \G' \stackrel{\tilde\xi'}{\to}{^L{\tilde{G}'}}\to{^L{\tilde{C}'}}\to 1,
\]
\item
every element of Out$_G(G')$ extends uniquely to an outer automorphism of $\tilde{G}'$ over $F$ which leaves $\tilde{C}'$ pointwise fixed.
\een
As a $K$-group, the $z$-extension $\tilde{G}'$ satisfies $\pi_0(\tilde{G}') = \pi_0(G')$, and $\tilde{G}'_\alpha$ is a $z$-extension of $G'_\alpha$ by $\tilde{C}'$ for each $\alpha\in \pi_0(G)$. Moreover, for any frame $(\psi',u')$ of $G'$ there is a corresponding frame $(\tilde\psi',\tilde{u}')$ for $\tilde{G}'$ such that $r_\alpha\tilde\psi_{\alpha\beta} = \psi_{\alpha\beta}r_\beta$ and $\tilde{u}_{\alpha\beta} = u_{\alpha\beta}$ for all $\alpha,\beta\in\pi_0(G)$.

Fix a central datum $(Z,\zeta)$, where $Z$ is a central induced torus of $G$, and $\zeta$ is a admissible character of $Z(F)$ if $F$ is local and a character of $Z(F)\backslash Z(\A)$ if $F$ is global. We also choose $\tilde Z$ and $\tilde\zeta$ on $\tilde G$ to be compatible with this datum.  Let $\tilde\eta'$ be the character dual to the Langlands parameter induced by the composition 
\[
W_F \to \G' \stackrel{\tilde\xi'}{\to}{^L{\tilde{G}'}}\to{^L{\tilde{C}'}},
\]
where $W_F\to \G'$ is any section. By condition (3), Out$_G(G')$ can be identified with a finite group of $F$-rational outer automorphisms of $\tilde{G}'$ which leave $\tilde{C}'$ pointwise invariant, thus fixing the central character $\tilde\eta'$.   We write $\tilde{Z}'$ for the preimage of $Z$ in $\tilde{G}'$, and $\tilde\zeta'$ for the product of $\tilde\eta'$ and the pullback of $\zeta$.  We can assume that the choice of auxiliary datum $(\tilde{G}',\tilde\xi')$ is compatible under isomorphisms of mesoscopic data $G'$, and therefore depends only on the elements $G'\in \F(G)$. 

Furthermore, we can also assume that $\tilde\xi'$ is of unitary type, in the sense that if $\phi':W_F\to \G'$ is an $L$-homomorphism such that the image of $\xi'\circ\phi'$ projects to a relatively compact subset of $\hat{G}$, then the image of $\tilde\xi'\circ\phi'$ also projects to a relatively compact subset of $\hat{\tilde{G}}$. The analogous condition ensures that the relative endoscopic transfer factors, defined for $K$-groups in \cite[\S2]{ArtTW}, have absolute value 1. 

\begin{rem}
In the endoscopic case, we have an endoscopic set $\E(G)$ consisting of ordinary endoscopic datum $(G^e,\G^e,s^e,\xi^e)$, where $s^e$ is a semisimple element in $\hat{G}$ satisfying certain assumptions. The associated auxiliary endoscopic datum $(\tilde{G}^e,\tilde\xi^e)$ is defined similarly, and also required to satisfy compatibility conditions, which we refer for example to \cite[\S2]{ArtTW} for details. In this paper, we shall generally indicate endoscopic objects with the superscript $^e$, and `stable' objects with $'$.
\end{rem}

\begin{rem}
We shall eventually show in Lemma \ref{indaux} that the stable transfer is independent of choice of auxiliary datum. On the other hand, recent work of Kaletha provides an alternative construction that amounts to a canonical choice of auxiliary datum \cite{covers}. 
\end{rem}

\section{Stable kernels and adjoint relations}
\label{gtf}

We now focus on the local setting. In preparation for the stable geometric transfer factors, we require several constructions related to the (inverse) Fourier transforms of stable orbital integrals and stable characters from \cite{LCR}, and develop some properties that we shall require. Most of this section is essentially review, aside from Lemma \ref{adj0}.

\subsection{Stable orbital integrals}
\label{soi}

Let $\C(G)$ be the space of Harish-Chandra Schwartz functions on $G(F)$, and let $\C(G,\zeta)$ be the subspace of $\zeta^{-1}$ equivariant functions, i.e., such that 
\[
f(xz) = \zeta(z)^{-1}f(x),\qquad x\in G(F), z\in Z(F).
\]
First, if $G$ is a connective reductive group, we define the normalised orbital integral
\[
f_G(\gamma) = |D(\gamma)|^{1/2}\int_{G_\gamma(F)\bs G(F)}f(x^{-1}\gamma x) dx, \qquad f\in\C(G,\zeta),
\]
where $D(\gamma)=\det(1-\text{Ad}(\gamma))_{\g/\g_\gamma}$ is the Weyl discriminant and $dx$ a fixed invariant measure on the orbit $G_\gamma(F)\bs G(F)$. If $G$ is a $K$-group, we set $f_G(\gamma) = f_{\alpha,G_\alpha}$ where $G_\alpha$ is the component that contains $\gamma$. Define 
\[
\I(G,\zeta) = \{f_G : f\in\C(G,\zeta)\},
\]
a topological space of functions on $\Gamma_\reg(G)$, topologised in a manner so that the map $f\to f_G$ is open and continuous. Denote by $\I_\cusp(G,\zeta)$ the subspace of cuspidal functions, that is, functions that vanish on the complement of $\Gamma_{\reg,\el}(G)$. We may view it as the space of functions that are annihilated by the restriction mapping $a_G\to a_M$ from $\I(G,\zeta)$ to $\I(M,\zeta)$ for any proper Levi $M$ of $G$. 
More precisely, let $P\in \P(M)$ such that $B\subset P$ and $T\subset M$ for a fixed pinning $(B,T,\{X_\alpha\})$ of $G$. Then $(B\cap M, T, \{X_{\alpha_M} \})$ is a pinning of $M$, where $\alpha_M$ runs over simple roots of $T$ relative to $M$. Fixing Haar measures on $G(F)$ and $M(F)$, we obtain measures on the unipotent radical $N_P$ and maximal compact subgroup $K$ by the formula
\[
\int_{G(F)} f(g) dg = \int_{M(F)} \int_{N_P(F)}\int_K f(muk) dk\ dn \ dm
\]
for all $f\in \C(G)$. The map $f\to f_M$ is then given by
\[
f_M(\gamma) = \int_{N_P(F)}\int_K f(k^{-1}n^{-1}\gamma nk) dn\ dk.
\] 
In particular, the map on the level of functions depends on the choice of $P$ and $K$, but choosing measures appropriately, it can be shown that the map induced from $\I(G,\zeta)$ to $\I(M,\zeta)$ indeed independent of these choices.

There is a natural measure on $\Gamma_\el(G)$ given by
\[
\int_{\Gamma_\el(G)}\alpha(\gamma)d\gamma = \sum_{\{T\}}|W(G(F),T(F))|^{-1}\int_{T(F)}\alpha(t)dt, 
\]
for any $\alpha\in C_c(\Gamma(G))$, where $\{T\}$ is a set of representatives of $G(F)$-conjugacy classes of elliptic maximal tori in $G(F)$, $W(G(F),T(F))$ is the Weyl group of $(G(F),T(F))$, and $dt$ is a fixed Haar measure on $T(F)$. The corresponding measures on $\Gamma_\el(M)$ determine a measure
\[
\int_{\Gamma(G)}\alpha(\gamma)d\gamma = \sum_{\{M\}}|W(M)|^{-1}\int_{\Gamma_\el(G)}\alpha(\gamma_M)d\gamma_M, 
\]
on $\Gamma(G)$.

The stable orbital integral of $f\in \C(G,\zeta)$ at $\delta\in\Delta_\reg(G)$ is given by
\[
f^G(\delta) = \sum_\gamma f^G(\gamma),
\]
where the sum is taken over the finite set of $\gamma\in\Gamma_\reg(G)$ that lie in the stable class $\delta$. We then define the subspace of $\I(G,\zeta)$
\[
S\I(G,\zeta)=\{f^G:f\in\C(G,\zeta)\},
\]
and set
\[
S\I_\cusp(G,\zeta) = \SI(G,\zeta) \cap \I_\cusp(G,\zeta).
\]
We call a tempered, $\zeta$-equivariant distribution on $G(F)$ stable if its value at any $f\in\C(G,\zeta)$ depends only on $f^G$. We similarly define measures on $\Delta_\el(G)$ and $\Delta(G)$ by
\[
\int_{\Delta_\el(G)}\beta(\delta)d\delta = \sum_{\{T\}_\text{st}}|W_F(G,T)|^{-1}\int_{T(F)}\beta(t)dt, 
\]
where $\beta\in C_c(\Delta(G)$, $\{T\}_\text{st}$ is a set of representatives of stable conjugacy classes of elliptic maximal tori in $G$ over $F$, and $W_F(G,T)$ is the subgroup of elements in the absolute Weyl group of $(G,T)$ defined over $F$; and
\[
\int_{\Delta(G)}\beta(\delta)d\delta = \sum_{\{M\}}|W(M)|^{-1}\int_{\Delta_\el(M)}\beta(\delta_M)d\delta_M.
\]
Note that for the induced torus $Z$, we have that $\Gamma(G)/Z(F) = \Gamma(G/Z)$ and $\Delta(G)/Z(F) = \Delta(G/Z)$. Also, we have $\Delta(\tilde G')/\tilde Z'(F) = \Delta(\tilde G'/\tilde Z') = \Delta(G').$

\subsection{Conjugacy classes again}
\label{conjcl}

 We shall construct certain `functorial' sets that keep track of the stable transfer mappings, parallel to $\Delta(G)$. Given $G$, the orbits of $\text{Out}_G(G')$ on $\Delta_{\el}(G')$ depend only on the isomorphism class of $G'$ in $\F(G)$. It makes sense then to define the set 
\[
\Delta^\F_{\el}(G) = \coprod_{G'\in \F_\el(G)}\Delta_{G,\el}(G')/\text{Out}_G(G'),
\]
which we can view as equivalence classes of pairs $(G',\delta')$. That is, if we write $\Delta_{G,\el}(G',G) = \Delta_{G,\el}(G')/\text{Out}_G(G')$, then 
\[
\Delta^\F_{\el}(G) = \{(G',\delta'): G'\in G'\in\F_\el(G), \delta' \in \Delta_{G,\el}(G')\}.
\]
We can similarly define $\Delta^\F_{G,\el}(M)$ for any $M\in\L$. Taking the union over the $W(M)$-orbits, we set
\[
\Delta^\F_{\reg}(G) = \coprod_{\{M\}}\Delta_{G,\el}^\F(M)/W(M).
\]
Fix an auxiliary datum $(\tilde{M'},\tilde\xi'_M)$ for $M'\in \F(M)$. We then define the set 
\[
\tilde\Delta_{G,\el}^\F(M) = \coprod_{M'\in\F_\el(M)}\Delta_{G,\el}(\tilde M'),
\]
which fibres over $\Delta_{G,\el}^\F(M)$, with the group $\prod_{M'}(\tilde{Z}'(F)\times \text{Out}_M(M'))$ acting transitively on the fibres. We again take the union of $W(M)$ orbits of $\tilde\Delta^\F_{G,\el}(M),$ and set
\[
\tilde\Delta^\F_{\reg}(G) = \coprod_{\{M\}}\tilde\Delta_{G,\el}^\F(M)/W(M),
\]
which fibres over $\Delta^\F_\reg(G)$. For brevity, we write $\Delta^\F_\reg(G)=\Delta^\F(G)$ and $\tilde\Delta^\F_\reg(G)=\tilde\Delta^\F(G)$. We can also view elements of $\tilde\Delta^\F(G)$ as equivalence classes of tuples $(G',\tilde{G}',\tilde\xi',\delta')$. 
These constructions are readily seen to generalise the `endoscopic' sets $\Delta^\E(G)$ and its variants in \cite[\S4]{STF1}, also denoted $\tilde{\Gamma}^\E(G)$ in \cite[\S2]{ArtTW}, which we shall use in this paper without comment.

\subsection{Endoscopic geometric transfer factors}
\label{egtf}
We briefly recall some basic facts about endoscopic transfer factors, such as in \cite[\S4-5]{STF1}. Given an endoscopic datum $G^e\in \E(G)$, the geometric endoscopic transfer factor is a smooth function $\Delta(\cdot,\cdot)$ on $\Delta_G(\tilde G^e)\times \Gamma(G)$ such as defined in \cite[\S2]{ArtTW}. The transfer factor determines a map
\[
f \to f^e(\delta^e) = \sum_{\gamma\in\Gamma(G)}\Delta(\delta^e,\gamma)f_G(\gamma),\qquad \delta^e\in\Delta_G(\tilde{G}^e)
\]
from functions $f\in \C(G,\zeta)$ to $f^e = f^{G^e}$ on $\Delta_G(\tilde{G}^e)$. The Langlands-Shelstad transfer then implies that $f^e$ belongs to $S\I(\tilde{G}^e,\tilde\zeta^e)$.  Fix an auxiliary endoscopic datum $(\tilde G^e,\tilde\xi^e)$ of $G^e$, so that $\tilde G^e$ is an extension of $G^e$ by a central induced torus $\tilde C^e$ with associated character $\tilde\eta^e$. The group $\tilde{C}^e(F)$ acts simply transitively on the fibres of the map $\Delta_G(\tilde{G}^e)\to \Delta_G(G^e)$, and $H^1(W_F,Z(\hat{\tilde{G}}^e))$ acts simply transitively on the set of $Z(\hat{\tilde{G}}^e)$-orbits of admissible embeddings $\tilde{\xi}^e$.  Then if $az\delta$ is the image in $\Delta^\E(G)$ of a point $(G^e,\tilde{G}^e,a\tilde\xi^e,z\delta^e)$ with $a \in H^1(W_F,Z(\hat{\tilde{G}}^e))$ and $z\in\tilde{C}^e(F)$, then the transfer factor satisfies 
$
\Delta(az\delta,\gamma)=\chi_a(\delta^e)\tilde\eta^e(z)\Delta(\delta,\gamma),
$
where $\chi_a$ is a character on $\tilde{G}^e$ determined by $a$ and the local Langlands correspondence for tori.

The transfer factors and consequently the Langlands-Shelstad transfer depend only on the image of $\delta^e$ in $\tilde\Delta^\E(G)$, and we can extend the transfer factors to $\tilde\Delta^\E(G)\times \Gamma(G)$ and define the extended map 
\[
f\to f^\E_G = \bigoplus_{G^e\in \E_\el(G)} f^e.
\]
That is, we define $\Delta(\delta^e,\gamma)$ to be zero unless there is an $M$ such that $(\delta^e,\gamma)$ belongs to the Cartesian product of $\tilde\Delta^\E_{G,\el}(M)/W(M)$ with $\Gamma_{G,\el}(M)/W(M)$. If there is such an $M$, then $(\delta^e,\gamma)$ is the image of a pair $(\delta^e_M,\gamma_M)$ in $\tilde\Delta^\E_{G,\el}(M)\times \Gamma_{G,\el}(M)$, and we set
\[
\Delta(\delta^e,\gamma) =\Delta_G(\delta^e,\gamma) = \sum_{w\in W(M)}\Delta_M(\delta^e_M,w\gamma_M).
\]
Each sum contains at most one nonzero term, and depends only on $\delta^e$ and $\gamma$.

Define also the adjoint transfer factor 
\[
\Delta(\gamma,\delta^e) = |\K_\gamma|^{-1}\overline{\Delta(\delta^e,\gamma)}
\]
on $\Gamma(G)\times \tilde\Delta^\E(G)$. Then by \cite[Lemma 2.3]{ArtTW}, we have the following adjoint relations
\[
\sum_{\delta^e\in\Delta^\E_\reg(G)}\Delta(\gamma,\delta^e)\Delta(\delta^e,\gamma_1) = \delta(\gamma,\gamma_1),\qquad \gamma,\gamma_1\in\Gamma(G),
\]
where $\delta(\cdot,\cdot)$ is the usual Kronecker delta, and
\[
\sum_{\gamma\in\Gamma_\reg(G)}\Delta(\delta^e,\gamma)\Delta(\gamma,\delta_1^e) =\tilde\delta(\delta^e,\delta_1^e),\qquad \delta^e,\delta_1^e\in\tilde\Delta^\E(G),
\]
where $\tilde\delta(\delta,\delta_1)=\tilde\eta^e(z)$ if $\delta_1=z\delta$ for some $z\in\tilde{C}^e(F)$ (or equivalently, if $\delta,\delta_1$ have the same projection onto $\Delta^\E_\reg(G)$) and equal to zero otherwise. The adjoint relations imply that $f_G\to f_G^\E$ is an isomorphism from $\I(G,\zeta)$ onto its image.

\subsection{Stable virtual characters}

Recall from \cite[\S3]{etc} the set $T(G)$ of $W_0$-orbits of essential triples $\tau = (L,\pi,r)$ where $L\in\L$, $\pi\in\Pi_2(L)$, and $r\in {R}_\pi$, where $\Pi_2(L)$ is the set of equivalence classes of irreducible unitary representations of $L(F)$ which are square integrable mod center, and ${R}_\pi$ is the $R$-group of $\pi$. Let $T_\el(G)$ be the subset of $\tau$ such that the kernel of $(1-r)$ acting on $\a_L$ is equal to $\a_G$. We define
\[
T(G) = \coprod_{\{M\}}T_\el(M)/W(M).
\]
We also have a decomposition with respect to any central induced torus $Z(F)$, which we assume contains the maximal $F$-split torus $A_G$,
\[
T(G) = \coprod_\zeta T(G,\zeta),
\]
where $\zeta$ runs over characters of $Z(F)$, and $ T(G,\zeta)$ is the subset of elements of $T(G)$ whose central character on $Z(F)$ equals $\zeta$. We also write $T_\el(G,\zeta) = T_\el(G)\cap T(G,\zeta)$.  The set $T(G)$ parametrises a family of locally integrable functions
\[
\gamma \to I(\tau,\gamma), \qquad \gamma\in\Gamma(G),
\]
such that for any $\zeta$, the functions $\overline{I(\tau,\gamma)}$ for $\tau\in T_\el(G,\zeta)$ form an orthogonal basis of $\I_\cusp(G,\zeta)$. Also, we have that $I(\tau,\gamma z) = I(\tau,\gamma)\zeta(z)$ for any $z\in Z(F)$. 

This family of functions has a stable analogue. If $F$ is nonarchimedean, by \cite[Lemma 5.1]{LCR} one can construct a set $\Phi_2(G,\zeta)$ parametrising a family of functions 
\[
\delta\to n(\delta)\overline{S(\phi,\delta)}, \quad \delta\in \Delta(G),
\]
which forms an orthogonal basis of $S\I_\cusp(G,\zeta)$, where $n(\delta)=|\K_\delta|$. Parallel to $T_\el(M)$ above, the basis then provides constructions of the larger sets 
\[
\Phi_2(G) = \coprod_\zeta \Phi_2(G,\zeta),\qquad \Phi(G) = \coprod_{\{M\}}\Phi_2(M)/W(M).
\]
The  set $\Phi(G)$ comes with an action $\phi\to\phi_\lambda=\phi\cdot\rho_\lambda$ where $\rho_\lambda$ in the nonarchimedean case is the unramified parameter which maps the Frobenius element to the image of $\lambda$ in $(Z(\hat{G})^\Gamma)^0$ under the exponential map, noting that $\a^*_{G,\CC}$ is equal to the Lie algebra of $(Z(\hat{G})^\Gamma)^0$. If $F$ is archimedean, we shall obtain this basis as a consequence of Lemma \ref{archkernel} below. 

In any case, if we take the local Langlands correspondence as known, we can identify  $\Phi_2(M,\zeta)$ with
the set of equivalence classes of cuspidal Langlands parameters $\phi:L_F\to {^LM}$ that are compatible with $\zeta$ in the sense that the composition of $\phi$ with the projection $^LM\to {^LZ}$ is the Langlands parameter defined by $\zeta$.  Here $L_F$ is $W_F$ or $W_F\times SL_2(\CC)$ depending on whether $F$ is real or $p$-adic. We can take instead
\[
S(\phi,\delta) = \sum_{\pi\in\Pi_2(G)}\Delta(\phi,\pi)I(\pi,\gamma),\qquad \gamma\in\Gamma(G)
\]
where $\Delta(\phi,\pi)$ are the endoscopic spectral transfer factors defined below, and $\Pi_2(G)$ is the set of equivalence classes of irreducible unitary representations of $G(F)$ that are square integrable mod centre. The orthogonality relation for $n(\delta)\overline{S(\phi,\delta)}$ defined this way follows from that of $I(\pi,\gamma)$.

The measure on $T_\el(G)$ is chosen to be
\be
\label{tell}
\int_{T_\el(G)}\alpha(\tau) d\tau = \sum_{\tau \in T_\el(G)/i\a^*_{G,\tau}}\int_{i\a^*_{G,\tau}}\alpha(\tau_\lambda)d\lambda
\ee
for any $\alpha\in C_c(T(G))$. Here we recall that $i\a^*_{G,F}=i\a^*_G/\a^\vee_{G,F}$, and also $i\a^*_{G,\tau} = i\a^*_G/\a^\vee_{G,\tau}$, where $\a^\vee_{G,\tau}$ is the stabiliser of $\tau$ in $i\a^*_G$, a lattice that lies between $\a^\vee_{G,F}$ and $\tilde\a^\vee_{G,F}$, and $d\lambda$ is a fixed measure on $i\a_{G,F}$. We then define the measure on $T(G)$ to be
\[
\int_{T(G)}\alpha(\tau)d\tau = \sum_{\{M\}}|W(M)|^{-1}\int_{T_\el(M)}\alpha(\tau_M)d\tau_M.
\]
Similarly, we define a measure on $\Phi_2(G)$ by setting
\be
\label{phi2}
\int_{\Phi_2(G)}\beta(\phi)d\phi = \sum_{\phi\in\Phi_2(G)/i\a^*_G}\int_{i\a^*_{G,\phi}}\beta(\phi_\lambda)d\lambda,
\ee
for any $\beta\in C_c(\Phi_2(G))$, where $i\a^*_{G,\phi} = i\a^*_G/\a^\vee_{G,\phi}$, where $\a^\vee_{G,\phi}$ is the stabiliser of $\phi$ in $i\a^*_G$. We then define the measure on $\Phi(G)$ to be
\be
\label{phim}
\int_{\Phi(G)}\beta(\phi)d\phi =  \sum_{\{M\}}|W(M)|^{-1}\int_{\Phi_2(M)}\beta(\phi_M)d\phi_M
\ee
for any $\beta\in C_c(\Phi_2(G))$.

\subsection{Endoscopic spectral transfer factors}

Now  for each elliptic endoscopic group $G^e\in \E_\el(G)$, we define 
\[
\Phi_2(\tilde{G}^e,G) = \Phi_2(\tilde{G}^e,\tilde\zeta^e)/\text{Out}_G(G^e).
\]
The spectral transfer factors $\Delta(\phi^e,\tau)$ are then defined in \cite[\S5]{STF1} to be uniquely determined functions on $\Phi_2(\tilde{G}^e,G)\times \tilde{T}_\el(G)$, satisfying 
\[
f^e(\phi^e)=\sum_{\tau\in T_\el(G)}\Delta(\phi^e,\tau)f_G(\tau),
\]
and $\Delta(\phi^e, z_\tau\tau)=\chi_\tau(z_\tau)\Delta(\phi^e,\tau)$ for $z_\tau\in Z_\tau$, where $Z_\tau=Z_\pi$ is a central subgroup used to define a central extension $\tilde R_\pi$ of $R_\pi$.  Define
\[
T^\E_\el(G) = \{(G^e,\phi^e): G^e\in \E_\el(G), \phi^e\in \Phi_2(\tilde G^e,G)\}
\]
and 
\[
T^\E(G) = \coprod_{\{M\}}T^\E_\el(M)/W(M).
\]
Then $\Delta(\phi^e,\tau)$ can be extended to a function on $T^\E(G)\times \tilde{T}(G)$ again as follows. We define $\Delta(\phi^e,\tau)$ to be zero unless there is an $M$ such that $(\phi^e,\tau)$ belongs to the Cartesian product of $T^\E_\el(M)/W(M)$ with $\tilde{T}_\el(M)/W(M)$. If there is such an $M$, then $(\phi^e,\tau)$ is the image of a pair $(\phi^e_M,\tau_M)$ in $T^\E_\el(M)\times \tilde{T}(M)$, and we set 
\[
\Delta(\phi^e,\tau)= \Delta_G(\phi^e,\tau) = \sum_{\tilde\tau_M}\Delta_M(\phi^e_M,\tilde\tau_M),
\]
where the sum runs over Weyl orbit $W(M)\tau_M$.

The orthogonality relation of the stable virtual characters which is a consequence of \cite[Lemma 5.1]{LCR} and its extension to strongly regular classes $\Delta(G)$ by \eqref{dell}, is given by 
\be
\label{ort1}
 \int_{\Delta_\el(G/Z)}n(\delta)S(\phi,\delta)\overline{S(\phi_1,\delta)}d\delta = \delta(\phi,\phi_1)n(\phi),
\ee
where $\delta$ is again the Kronecker delta, and $n(\phi)$ is simply defined to be the value 
\[
n(\phi)= \int_{\Delta_\el(G/Z)}n(\delta)S(\phi,\delta)\overline{S(\phi,\delta)}d\delta .
\]
It is parallel to the formula for the invariant virtual characters
\[
\int_{\Gamma_\el(G/Z)}I(\tau,\gamma)\overline{I(\tau_1,\gamma)}d\gamma = \delta(\tau,\tau_1)n(\tau),
\]
where the constant $n(\tau)$ is defined by \cite[Theorem 6.2]{etc}. We shall later derive a dual orthogonality relation for $S(\delta,\phi)$ in Lemma \ref{ort2} below.

\subsection{Fourier transforms}

To define our stable transfer factors, we must recall some constructions relating to the (inverse) Fourier transforms of orbital integrals.  We have for any $f\in \H(G,\zeta)$, the relations
\be
\label{fggamma}
f_G(\gamma) = \int_{T(G,\zeta)}I(\gamma,\tau)f_G(\tau)d\tau
\ee
and
\[
f_G(\tau) = \int_{\Gamma(G/Z)}I(\tau,\gamma)f_G(\gamma)d\gamma,
\]
where we denote by $I(\tau,\gamma)= |D(\gamma)|^\frac12\Theta(\tau,\gamma)$ the normalised virtual character associated to $\tau$, and $I(\gamma,\tau)$ on the other hand can be viewed as the coefficient in the Fourier inversion of the orbital integral $f_G(\gamma)$. They are smooth functions in both variables, described in Theorems 4.1 and 4.3 of \cite{fourier}. We shall be interested in their stable analogues. By \cite[Lemma 6.3]{LCR} there exist smooth functions $S(\delta,\phi)$ and $S(\phi,\delta)$ of $\phi\in\Phi(G,\zeta)$ and $\delta\in\Delta(G/Z)$, which are respectively $\zeta$ and $\zeta^{-1}$-equivariant under translation by $Z(F)$, such that
\be
\label{fgdelta}
f^G(\delta) = \int_{\Phi(G,\zeta)}S(\delta,\phi)f^G(\phi)d\phi
\ee
and
\be
\label{fgphi}
f^G(\phi) = \int_{\Delta(G/Z)}S(\phi,\delta)f^G(\delta)d\delta
\ee
for any $f\in \H(G,\zeta)$. The smooth functions are given by
\be
\label{Sformula1}
S(\delta,\phi) = \sum_{\gamma\in\Gamma(G)}\sum_{\tau\in T(G,\zeta)}\Delta(\delta,\gamma)I(\gamma,\tau)\Delta(\tau,\phi)
\ee
and 
\be
\label{Sformula2}
S(\phi,\delta) = \sum_{\tau \in T(G,\zeta)}\sum_{\gamma\in\Gamma(G)}\Delta(\phi,\tau)I(\tau,\gamma)\Delta(\gamma,\delta).
\ee
Here $\Delta(\delta,\gamma)$ is the endoscopic geometric transfer factor with adjoint $\Delta(\gamma,\delta)$, and $\Delta(\phi,\tau)$ is the endoscopic spectral transfer factor with adjoint $\Delta(\tau,\phi)$, as recalled above. While it is probably best to renormalise these transfer factors according to the works of Kaletha (c.f. \cite[\S4]{kal}), we neglect to to do so here.

We record here the archimedean analogue.

\begin{lem}
\label{archkernel}
Let $F$ be an archimedean local field. Then there exist smooth functions $S(\phi,\delta)$ and $S(\delta,\phi)$ of $\phi\in\Phi(G,\zeta)$ and $\delta\in\Delta(G)$, which are respectively $\zeta$ and $\zeta^{-1}$-equivariant under translation of $\delta$ by $Z(F)$, such that 
\[
f^G(\delta) = \int_{\Phi(G,\zeta)}S(\delta,\phi)f^G(\phi)d\phi, \qquad f^G(\phi) = \int_{\Delta(G/Z)}S(\phi,\delta)f^G(\delta)d\delta
\]
for any $f\in \H(G,\zeta)$. 
\end{lem}

\begin{proof}
As in the nonarchimedean case, the proof will follow in the same way as \cite[Lemma 6.3]{LCR} from the property that the linear mapping 
\[
f\to f^G_\text{gr}(\phi) = \sum_{\tau\in T(G,\zeta)}\Delta(\phi,\tau)f_G(\tau),\qquad \phi\in\Phi(G,\zeta)
\]
is stable, and induces an isomorphism from  $S\I(G,\zeta)$ to the graded vector space
\[
S\I_\text{gr}(G,\zeta) = \bigoplus_{\{M\}}\SI_\cusp(M,\zeta)^{W(M)},
\]
(see also \eqref{sigr} below). These in turn can be deduced from \cite[IV.2]{MW1}, in particular Th\'eor\`eme IV.2.3(i) and Corollaire IV.2.9. Moreover, the functions once again take the form \eqref{Sformula1} and \eqref{Sformula2}.
\end{proof}

\subsection{An adjoint relation for stable kernels}
\label{adjapp}

We derive the following identity which will play an important role in proving identities for our stable transfer factor. It is the stable analogue of the relation \eqref{tau1} in \cite[Theorem 4.5]{fourier} relating $I(\tau,\gamma)$ to $I(\gamma,\tau)$ for elliptic virtual characters.

\begin{lem}
\label{adj0}
The stable kernels satisfy the adjoint relation
\be
\label{adjoint}
n(\delta) \overline{S(\phi,\delta)} = n(\phi) S(\delta,\phi)
\ee
for $\phi\in\Phi_2(G,\zeta)$ and $\delta\in \Delta_\el(G)$.
\end{lem}

\noindent 
We can motivate the identity as follows. Applying the inversion formulae \eqref{fgphi} and \eqref{fgdelta} consecutively to 
\[ 
f^G(\phi) = \int_{\Delta(G/Z)}S(\phi,\delta)\int_{\Phi(G,\zeta)}S(\delta,\phi_1)f^G(\phi_1)d\phi_1 d\delta,
\]
and then interchanging the integrals, we have
\[
\int_{\Phi(G,\zeta)}\int_{\Delta(G/Z)}S(\phi,\delta)S(\delta,\phi_1)d\delta\ f^G(\phi_1)d\phi_1 .
\]
Then if we had the relation \eqref{adjoint}, this is 
\[
\int_{\Phi(G,\zeta)}n(\phi)^{-1}\int_{\Delta(G/Z)}n(\delta)S(\phi,\delta)\overline{S(\phi_1,\delta)}d\delta\ f^G(\phi_1)d\phi_1,
\]
so that the orthogonality relation \eqref{ort1} gives us the tautaology
\[
\int_{\Phi(G,\zeta)}\delta(\phi,\phi_1)f^G(\phi_1)d\phi_1 = f^G(\phi),
\]
as expected.

\begin{proof}
For simplicity, we assume $(Z,\zeta)$ to be trivial, since it does not affect the proof. We would like to compare
\[
S(\delta,\phi) = \sum_{\gamma\in\Gamma(G)}\sum_{\tau\in T(G,\zeta)}\Delta(\delta,\gamma)I(\gamma,\tau)\Delta(\tau,\phi)
\]
with
\[
\overline{S(\phi,\delta)} = \sum_{\tau \in T(G,\zeta)}\sum_{\gamma\in\Gamma(G)}\overline{\Delta(\phi,\tau)I(\tau,\gamma)\Delta(\gamma,\delta)}.
\]
We recall the identities satisfied by the endoscopic geometric and spectral transfer factors, relating them to their adjoint functions from equations (2.3) and (5.5) in \cite{LCR},
\be
\label{rel1}
\Delta(\delta^e,\gamma) = n(\delta)\overline{\Delta(\gamma,\delta^e)} 
\ee
and
\be
\label{rel2}
\Delta(\tau,\phi^e) = |Z(\hat{G}^e)^\Gamma/Z(\hat{G})^\Gamma|^{-1}n(\tau)n(\phi)^{-1}\overline{\Delta(\phi^e,\tau)},
\ee
where 
\[
n(\tau) = |R_{\pi,r}||\det(1-r)_{\a_L/\a_G}|,
\]
and $R_{\pi,r}$ is the centraliser of $r$ in the $R$-group $R_\pi$. Recall that $T(G)$ is the set of $W_0^G$-orbits of essential triplets $\tau = (L,\pi,r)$ where $L\in\L$, $\pi$ is an equivalence class of an irreducible unitary representation of $L(F)$ that is square integrable modulo centre, and $r\in R_\pi$. The $R$-group of $\pi$ is the quotient $R_\pi = W_\pi/W^0_\pi$, where $W_\pi^0$ is the subgroup of elements of $w\in W_\pi$ such that the normalised intertwining operator $R(\pi,w)$ acts by a scalar.  The subset $T_\el(G)$ consists of $\tau$ for which the kernel of $(1-r)$ acting on $\a_L$ equals $\a_G$. Finally, we note that in the case at hand, we shall identify $\delta$ with its image $\delta^e = \delta^*$ in $\Delta_G(G^*)$, similarly $\phi$ with $\phi^e=\phi^*$ in $\Phi(G^*,\zeta^*)$.

Moreover, if $\tau\in T_\el(G)$ we write $\tau^\vee = (L,\pi^\vee,r)$ for the contragredient and set 
\[
i^G(\tau) = |\det(1-r)_{\a_L/\a_G}|^{-1},
\]
then it follows from the special case $M=G$ of \cite[Theorem 4.5]{fourier} that
\be
\label{tau1}
I^\text{old}(\gamma,\tau) = i^G(\tau) I(\tau^\vee,\gamma).
\ee
for any $\gamma \in \Gamma_\el(G)$ and $\tau\in T_\disc(G)$, which we define below. But it is crucial that the measures on $T(G)$ assigned in \cite[\S4]{fourier} differs from that of \cite[\S4]{LCR}  by a factor of $|R_{\pi,r}|$, which we must reconcile. This explains our notation $I^\text{old}(\gamma,\tau)$ for the kernel used in \cite{fourier}. We first explain the measures. We write $T_\disc(G)$ for the subset of $W_0^G$-orbits for which the set of regular elements 
\[
W_\pi(r)_\reg = \{w \in W_\pi(r) : \a^w_{L}= \a_G\}
\]
is nonempty. Here $W_\pi(r)$ is the subset of elements in $W(\a_L) = W^0_\pi\cdot r$ which stabilise $\pi$ and which have the same projection onto the $R$-group as $r$. For any $w$ in this set, we write $\varepsilon_\pi(w)$ for the sign of the element $wr^{-1}$ in the Weyl group $W^0_\pi$. The function $i(\tau)= i^G(\tau)$ is more generally defined on $T_\disc(G)$ as
\[
i(\tau) = |W^0_\pi|^{-1}\sum_{w \in W_\pi(r)_\reg }\varepsilon_\pi(w)|\det(1-w)_{\a_L/\a_G}|^{-1}. 
\]
For $\tau\in T_\el(G)$, the group $W^0_\pi$ is trivial and $i(\tau)$ specialises to the former expression. Now the measure on $T_\disc(G)$ is chosen in \cite{fourier} to be 
\[
\int_{T_\disc(G)}\alpha(\tau)d\tau = \sum_{\tau \in T_\disc(G)/i\a^*_G}|R_{\pi,r}|^{-1}|\a^\vee_{G,\tau}/\a^\vee_{G,F}|^{-1}\int_{i\a^*_{G,F}}\alpha(\tau_\lambda)d\lambda
\]
for any $\alpha\in C_c(T_\disc(G))$. On the other hand, recalling the measure chosen in \eqref{tell} compatibly with \cite{LCR}, it follows that the righthand side can then be written as the sum over $\tau \in T_\el(G)/i\a^*_{G}$ of 
\[
|i\a^*_G/\a^\vee_{G,\tau}|^{-1}\int_{i\a^*_{G,\tau}}\alpha(\tau_\lambda)d\lambda = |\a^\vee_{G,\tau}/\a^\vee_{G,F}|^{-1}\int_{i\a^*_{G,F}}\alpha(\tau_\lambda)d\lambda.
\]
In particular, we see that the measure conversion from \cite{fourier} to \cite{LCR} is given by multiplication by $|R_{\pi,r}|^{-1}$. 

We shall show that the factor $i(\tau)$ in \eqref{tau1} must therefore be multiplied by the same factor in order for the choice of measures to be consistent. We shall in fact prove a slightly stronger result, that is, for the weighted kernels $I_M(\gamma,\tau)$, by which $I_G(\gamma,\tau)= I(\gamma,\tau)$ is a special case. By our choice of measure, \cite[Theorem 4.1]{fourier} asserts the existence of a smooth function $I_M(\gamma,\tau)$ on $\gamma\in \Gamma(M)\cap G_\reg(F)$ and $\tau\in T_\disc(L)$ for $L\in \L$ such that
\[
I_M(\gamma,f) = \sum_{L\in \L}|W^L_0||W^G_0|^{-1}\int_{T_\disc(L)}|R_{\pi,r}|^{-1}I_M^\text{old}(\gamma,\tau)f_L(\tau)d\tau,
\]
where $I_M(\gamma,f)$ is the weighted orbital integral of $f\in \C(G)$, and $R_{\pi,r}$ is the group associated to $\tau=(L,\pi,r)$. In particular, the kernel $I_M^\text{old}(\gamma,\tau)$ of \cite{fourier} relates to the kernel $I^M(\gamma,\tau)= I_M^\text{new} (\gamma,\tau)$ of \cite{LCR}, which is the one we are using, by the renormalisation
\[
I^M(\gamma,\tau) = |R_{\pi,r}|^{-1}I_M^\text{old}(\gamma,\tau).
\]
We may as well verify the formula by recalling Arthur's argument. Substituting this expression into the geometric side of the invariant local trace formula \cite[(5.1)]{fourier},
\[
\sum_{M\in \L}|W^M_0||W^G_0|^{-1}(-1)^{\dim(A_M/A_G)}\int_{\Gamma_\el(M/Z)}I_M(\gamma,f) g_M(\gamma)d\gamma
\]
for $g\in C_c^\infty(G_\reg(F))$, we obtain the expression
\[
\sum_{L\in \L}|W^L_0||W^G_0|^{-1}(-1)^{\dim(A_L/A_G)}\int_{T_\disc(L)}|R_{\pi,r}|^{-1}I_M'(\tau,g)f_L(\tau^\vee)d\tau,
\]
where
\[
I'_L(\tau,g) = \sum_{M\in \L}|W^M_0||W^G_0|^{-1}(-1)^{\dim(A_M\times A_L)}\int_{\Gamma_\el(M/Z)}I_M^\text{old}(\gamma,\tau^\vee) g_M(\gamma)d\gamma.
\]
By the local trace formula, this is equal to the spectral expansion
\[
\sum_{L\in \L}|W^L_0||W^G_0|^{-1}(-1)^{\dim(A_L/A_G)}\int_{T_\disc(L)}|R_{\pi,r}|^{-1}i^L(\tau)I_L(\tau,g)f_L(\tau^\vee)d\tau.
\]
The remainder of the argument follows that of \cite[\S6]{fourier}. Namely, considering the difference of the spectral and geometric expansions as distributions in $f_G$, we see that the difference is a finite sum of smooth symmetric functions on the strata $T_\disc(L)$ of $T(G)$ as $L$ varies. Since $f_G$ ranges over $\I(G)$, we can separate the contributions of the various strata, and it follows that 
\[
i^L(\tau) I_L(\tau,g) = I'_L(\tau,g),\qquad L\in \L, \tau\in T_\disc(L).
\]
From this, the same argument as Arthur's gives the parallel expansion
\[
I_L(\tau,g) =  \sum_{M\in \L}|W^M_0||W^G_0|^{-1}\int_{\Gamma_\el(M/Z)}I_L(\tau,\gamma)g_M(\gamma)d\gamma
\]
of \cite[(4.1)$^\vee$]{fourier}. Then comparing the expansions for $I_L(\tau,g)$ and $ I'_L(\tau,g)$ that we have obtained, we have again that  
\[
(-1)^{\dim(A_M\times A_L)}I_M^\text{old}(\gamma,\tau^\vee) = i^L(\tau)I_L(\tau,\gamma).
\]
Using the fact that $i^L(\tau^\vee) = i^L(\tau)$ and setting $M=L=G$, we conclude that the identity \eqref{tau1} should be indeed multiplied by $|R_{\pi,r}|^{-1}$ to give
\be
\label{igt}
I(\gamma,\tau) = |R_{\pi,r}|^{-1}i^G(\tau) I(\tau^\vee,\gamma)
\ee
by our choice of measures.

Now we can prove the proposition. Combining the identities together \eqref{rel1}, \eqref{rel2}, and \eqref{igt}, it follows that $S(\delta,\phi)$ is equal to
\[
\frac{n(\delta)}{|Z(\hat{G}^e)^\Gamma/Z(\hat{G})^\Gamma|n(\phi)}\sum_{\tau \in T(G)}\sum_{\gamma\in\Gamma(G)} n(\tau)|R_{\pi,r}|^{-1}i^G(\tau) \overline{\Delta(\gamma,\delta) }{I(\tau^\vee,\gamma)}\overline{\Delta(\phi,\tau)},
\]
which simplifies to
\[
{n(\delta)}{n(\phi)^{-1}}\sum_{\tau \in T(G)}\sum_{\gamma\in\Gamma(G)}  \overline{\Delta(\gamma,\delta) }{I(\tau^\vee,\gamma)}\overline{\Delta(\phi,\tau)},
\]
where we have used the fact that the quotient $Z(\hat{G}^e)^\Gamma/Z(\hat{G})^\Gamma$ is trivial for $G^e = G^*$.  Finally, we see that if 
\[
I(\tau^\vee,\gamma)  = \overline{I(\tau,\gamma)},
\]
then the desired formula follows. We simply deduce this from the properties that
\[
\overline{f_G(\gamma)} = \overline{f}_G(\gamma),\qquad \overline{f_G(\tau^\vee)} = \bar{f}_G(\tau),
\]
(see for example, the proof of \cite[Theorem 6.1]{etc}) and comparing the expansions on the either side of the first identity using \eqref{fggamma},
\[
\int_{T(G)}\overline{I(\gamma,\tau)f_G(\tau)}d\tau = \int_{T(G)}I(\gamma,\tau)\overline{f}_G(\tau)d\tau = \int_{T(G)}I(\gamma,\tau^\vee)\overline{f_G(\tau)}d\tau,
\]
and again varying $f_G$ in $\I(G)$ accordingly. From this we have that 
\[
I(\gamma,\tau^\vee)  = \overline{I(\gamma,\tau)},
\]
then using the fact that $i_G(\tau^\vee)= i_G(\tau)$ and the relation \eqref{igt}, the claim follows.

Note that from the proof above we also find the parallel statement for the invariant kernels, namely, $I(\gamma,\tau) = \tilde{i}^G(\tau) \overline{I(\tau,\gamma)}$ where $\tilde{i}^G(\tau)=|R_{\pi,r}|^{-1}i^G(\tau)$.  

\end{proof}
At various times in the following sections, we will want to interchange the order of integration of stable kernels in a manner similar to the heuristic above. We will give a justification for it in \S\ref{proof}.

\section{Stable geometric transfer factors}
\label{sgtf}

\subsection{Local Langlands correspondence}
Our construction of stable transfer factors relies on the transfer of $L$-parameters. Thus it is necessary to assume  the local Langlands conjecture in order to formulate that definition. The local Langlands group $L_F$ is defined to be $W_F$ if $F$ is archimedean and $W_F\times SU(2)$ if $F$ is nonarchimedean. Recall that an $L$-homomorphism in this context is a homomorphism $\phi:L_F\to {^LG}$ that commutes with projections onto $W_F$ of its source and target. We say it is admissible if it is continuous and sends elements of $W_F$ to semisimple elements of $^LG$, and relevant if its image being contained in a Levi subgroup $^LM$ of $^LG$ implies that $^LM$ is the $L$-group of a Levi subgroup $M$ of $G$ over $F$.

Let $\Phi^+(G)$ be the set of $\hat{G}$-conjugacy classes of relevant, admissible $L$-homomorphisms $\phi$.  By abuse of notation, let 
$
\Phi(G) = \Phi_\temp(G)
$
be the subset of bounded or tempered Langlands parameters $\phi\in\Phi(G)$, that is, whose image projects onto a relatively compact subset of $\hat{G}$. Their corresponding $L$-packets $\Pi_\phi$ are expected to consist of tempered representations. We also let $\Phi_2(G)$ be the set of cuspidal parameters $\phi$, whose image does not lie in a proper parabolic subgroup $^LP$ of $^LG$. (The sets described above by the same notation were introduced by Arthur in order to avoid the use of the local Langlands correspondence.) The local Langlands conjecture, at the most basic level, asserts that the set $\Pi^+(G)$ of irreducible admissible representations of $G(F)$ can be written as a disjoint union of finite packets $\Pi_\phi$ as $\phi$ ranges over $\Phi^+(G)$. In other words, there exists a surjective map 
$
\Pi(G)^+\to \Phi^+(G)
$
with finite fibres, which restricts to a surjective map of tempered representations to tempered parameters $\Pi(G)=\Pi_\temp(G)\to \Phi_\temp(G)$.  

In the present context, given $G' \in \F(G)$ with auxiliary datum $(\tilde{G}',\tilde\xi')$, the associated Langlands parameters are $\hat{G}$-conjugacy classes of homomorphisms $W_F \to \G' $, and composing with $\tilde\xi'$ we have a map into an $L$-group $^L{\tilde G}'$. We shall write $\Pi(G,\zeta)$ for the subset of representations with central character equal to $\zeta$, and similarly $\Phi(G,\zeta)$, whereby
\[
\Phi(G) = \coprod_\zeta \Phi(G,\zeta).
\]
By construction, for any $G'$, the auxiliary data $(\tilde{G}',\tilde\xi')$ and central datum $(Z,\zeta)$ are chosen such that any $\phi'\in\Phi(\tilde{G}',\tilde\zeta')$ maps to a parameter $\phi\in\Phi(G,\zeta)$, again with parallel restrictions to tempered parameters.

The existence of this surjective map for $\Phi(G)$ (and the releveant endoscopic character identities) allows us to deduce the existence stable transfer mapping. In particular, as a result of the refined local Langlands correspondence, the set $\Phi(G)$ parametrises tempered $L$-packets of $G(F)$, and the space of stable orbital integrals on regular semisimple elements of $G(F)$ corresponds to the Paley-Wiener space on $\Phi(G)$ under the map given by taking stable characters. 

The local Langlands correspondence implies stable transfer. The proof follows \cite[p.999]{Mok} and relies on the trace Paley-Wiener theorem for Schwartz functions \cite{PW}. This result was recently established for $G'$ of general rank \cite{sunohara}, under the nonarchmedean local Langlands correspondence. Nonetheless, we include the proof below for the equal rank case to illustrate the main idea.

\begin{prop}
\label{thm1}
Assume the refined local Langlands correspondence for $G$ and $G'$ over $F$ nonarchimedean. Then for any local field $F$ of characteristic zero and $f\in \C(G,\zeta)$, there exists a unique $f'\in \SI(\tilde G',\tilde\zeta')$ characterised by \eqref{bemap}.
\end{prop}

\begin{proof}
If $F$ is nonarchimedean, the action $\phi\to \phi_\lambda$ of $i\a^*_G$ on $\Phi_2(G)$ makes it into a disjoint union of compact tori of the form $i\a^*_{G,\phi} = i\a^*_G/\a^\vee_{G,\phi}$, where $\a^\vee_{G,\phi}$ is the stabiliser of $\phi$ in $i\a^*_G$. The orthogonal basis $n(\delta)\overline{S(\phi,\delta)}$ makes $\SI_\cusp(G,\zeta)$ into the Paley-Wiener space on $\Phi_2(G,\zeta)$, in the sense that it is the space of functions on $\Phi_2(G,\zeta)$ supported on finitely many connected components, and which on the component of any $\phi$ pullback to a finite Fourier series on $i\a^*_{G,\phi}$ \cite[p.541]{LCR}. Moreover, the larger graded vector space
\be
\label{sigr}
\SI_\text{gr}(G,\zeta) = \bigoplus_{\{M\}}\SI_\cusp(M,\zeta)^{W(M)}
\ee
can be identified with the natural Paley-Wiener space on $\Phi(G,\zeta)$. It is a consequence of \cite[Theorem 6.1]{LCR} that we may identify $\SI_\text{gr}(G,\zeta)$ with $\SI(G,\zeta)$. If $F$ is archimedean, we recall that $\Phi(G,\zeta)$ is the space of $\zeta$-equivariant tempered Langlands parameters, and is a basis for $\SI(G,\zeta)$.

Now for general $F$, the space $\SI(G,\zeta)$ corresponds to the Paley-Wiener space of $\Phi(G,\zeta)$ by the map 
\[
\phi \to f^G(\phi), \qquad f\in\C(G,\zeta).
\]
Moreover, the $L$-embedding $\tilde\xi'$ induces isomorphisms on the corresponding maximal tori, and it follows that the function 
\[
\phi'\to f^G(\tilde\xi'\circ\phi'),\qquad \phi'\in \Phi(\tilde G')
\] 
belongs to the Paley-Wiener space on the set of tempered local Langlands parameters of $\tilde G'(F)$. Then there exists a function $f' \in \C(\tilde G',\tilde\zeta')$, uniquely determined up to its stable orbital integral, such that $f'(\phi') = f^G(
\tilde\xi'\circ\phi')$.
\end{proof}

\begin{rem}
Suppose $F$ is a nonarchimedean field over which $G$ is unramified. Let $\H(G,K)$ denote the spherical Hecke algebra of $G(F)$, where $K$ is a hyperspecial maximal compact subgroup of $G(F)$. The Satake isomorphism implies an algebra homomorphism of Hecke algebras induced by the restriction of the embedding $\xi':\G'\to {^LG}$, compatible with this spectral mapping $f\to f'$. 
\end{rem}

The goal of the stable transfer factors is to provide a parallel construction in terms of stable orbital integrals. We first note that the stable spectral transfer immediately implies an abstract stable geometric transfer.

\begin{cor}
\label{skt}
Given \eqref{bemap}, there exists a stable distribution $\Theta$ on $G\times G'$ such that
 \[
 f'(\delta') = \int_{\Delta(G/Z)} \Theta(\delta',\delta)f^G(\delta) d\delta.
 \]
\end{cor}

\begin{proof}
The stable transfer mapping of Proposition \ref{thm1} gives a continuous linear map from $\C(G,\zeta) \to \SI(G',\zeta')$. 
In particular, the $K$-invariant functions in the Harish-Chandra Schwartz space form a nuclear Fr\'echet space, so that $\C(G,\zeta)$ is an LF space, that is, a countable strict inductive limit of Fr\'echet spaces, in the case where $F$ is nonarchimedean. If $F$ is archimedean, it is simply nuclear Fr\'echet \cite[\S51]{treves}.  Then an application of the Schwartz kernel theorem \cite[Ch. II, \S3, Th\'eor\`eme 12]{groth} (also \cite[Theorem 51.7]{treves}) applied to nuclear LF spaces allows us to identify this mapping with the existence of the integral kernel as desired.
\end{proof}

One main goal of this paper is to propose an explicit construction for this kernel. 

\subsection{Stable geometric transfer factors}

We write $S'(\cdot,\cdot) = S^{\tilde{G}'}(\cdot,\cdot)$ for the kernel functions associated to ${G}' \in \F(G)$.  We can now introduce our stable geometric transfer factor as a stable distribution on $\Delta_G(\tilde{G}')\times \Delta(G)$,
\begin{equation}
\label{candidate}
\Theta_{\tilde\xi'}(\delta',\delta) =\int_{\Phi(\tilde G',\tilde\zeta')}S'(\delta',\phi')S(\tilde\xi'\circ\phi',\delta)d\phi',
\end{equation}
where we identify $\tilde\xi'\circ\phi'$ with the image of $\phi'$ in $\Phi(G,\zeta)$ determined by $\xi'$. (Recall that the stable character $f^G(\phi)$ is independent of choice of auxiliary datum.) We distinguish this distribution from the abstract one in Corollary \ref{skt} by the subscript. Notice that unlike the endoscopic transfer factor, we cannot require that $\Theta_{\tilde\xi'}$ vanish if $\delta'$ is not an image (we are indebted to the referee for this observation). By construction, it is clear that the stable transfer factor depends only on the stable conjugacy classes of its inputs. It also implicitly depends on the normalisation of the Langlands-Shelstad transfer factors.  

The stable transfer factors can be extended to distributions on $\tilde\Delta^\F(G)\times \Delta(G)$ as follows. Set $\Theta_{\tilde\xi'}(\delta',\delta)$ to be zero unless there is an $M$ such that $(\delta',\delta)$ belongs to the Cartesian product of $\tilde\Delta^\F_{G,\el}(M)/W(M)$ with $\Delta_{G,\el}(M)/W(M)$. If there is such an $M$, then $(\delta',\delta)$ is the image of a pair $(\delta'_M,\delta_M)$ in $\tilde\Delta^\F_{G,\el}(M)\times \Delta_{G,\el}(M)$, and we set
\be
\label{extended}
\Theta_{\tilde\xi'}(\delta',\delta) = \Theta_{\tilde\xi',G}(\delta',\delta) =\sum_{w\in W(M)}\Theta_{\tilde\xi',M}(\delta'_M,w\delta_M),
\ee
where again each sum contains at most one nonzero term, and depends only on $\delta'$ and $\delta$. 

The candidate kernel is the main construction of this paper, which we conjecture realizes the kernel of stable geometric transfer, and which we attribute to Shelstad \cite{shelstad}. Note that as in the case of endoscopy, the stable transfer factor is not unique as it depends on various choices of data that we shall elaborate in the next section. 

\begin{con}
\label{maincon}
The distribution \eqref{candidate} is a kernel satisfying Corollary \eqref{skt}. That is, we may take $\Theta= \Theta_{\tilde\xi'}$.
\end{con}

We shall use our stable transfer factor \eqref{candidate} to study the stable transfer of orbital integrals \eqref{bemap2}, and relate it to the stable transfer of characters \eqref{bemap}. But first, we shall first develop some basic properties that are parallel to that of the endoscopic transfer factors.

\begin{rem}
In \cite{thomas,shelstad} it is explicitly computed that the stable transfer factor for $G = \text{SL}_2(\R)$ and $G' = \text{SO}_2(\R)$ gives a divergent infinite sum, but can be used to produce the stable transfer of compactly-supported functions.
\end{rem}


\begin{lem}\label{lem1}
Let $G'\in \F_\el(G)$. 
\begin{enumerate}
\item[(i)]
For any $z\in\tilde{Z}'(F)$ whose image in $Z(F)$ equals $z_G$, we have
\[
\Theta_{\tilde\xi'}(\delta' z ,\delta z_G ) = \tilde\zeta'(z)^{-1}\Theta_{\tilde\xi'}(\delta',\delta)\zeta(z_G).
\]
\item[(ii)]
Given the injective linear map $\lambda\to \lambda'$ from $\a^*_{G,\CC}$ to $\a^*_{\tilde{G}',\CC}$, we have 
\[
e^{\lambda'(H_{\tilde{G}'}(\delta'))}\Theta_{\tilde\xi'}(\delta',\delta) = e^{\lambda(H_{G}(\delta))}\Theta_{\tilde\xi'}(\delta',\delta).
\]
whenever $\delta'$ is an image of $\delta$.
\end{enumerate}
\end{lem}

\begin{proof}
From the definition, we first write
\[
\Theta_{\tilde\xi'}(\delta' z ,\delta z_G ) =\int_{\Phi(\tilde G',\tilde\zeta')}S'(\delta' z,\phi')S(\tilde\xi'\circ\phi',\delta z_G)d\phi',
\]
and the equivariance properties of the stable kernels (which follow from that of the endoscopic transfer factors and invariant kernels, c.f. \cite[Lemma 6.3]{LCR}) yield 
\[
S'(\delta' z,\phi')S(\tilde\xi'\circ\phi',\delta z_G) = \tilde\zeta'(z)^{-1}S'(\delta' ,\phi')S(\tilde\xi'\circ\phi',\delta)\zeta(z_G).
\]
The first result follows.

In the second place, the map from $\a^*_{G,\CC}$ to $\a^*_{\tilde{G}',\CC}$ can be viewed as a map of the complex Lie algebras of 
\[
(Z(\hat{G})^\Gamma)^0\stackrel{\sim}{\to}(Z(\hat{\tilde{G}}')^\Gamma)^0,
\]
where the isomorphism follows from the fact that $G'$ is elliptic. It follows then that there is an injection from  $(Z(\hat{G}')^\Gamma)^0$ to  $(Z(\hat{\tilde{G}}')^\Gamma)^0$ dual to the projection $\tilde{G}'\to G$ given by property (1) of the auxiliary datum. Then the second identity follows.
\end{proof}

\subsection{Adjoint relations}

Define the adjoint stable geometric  transfer factor 
\be
\label{adjdef}
\Theta_{\tilde\xi'}^\mathrm{ad}(\delta',\delta) = n(\delta')^{-2}\overline{\Theta_{\tilde\xi'}(\delta',\delta)},
\ee
and its extension analogous to \eqref{extended}.The stable geometric transfer factors then satisfy adjoint relations parallel to those of the endoscopic geometric transfer factors in Section \ref{egtf}. But first, we derive an orthogonality relation for the $S(\delta,\phi)$ analogous to \eqref{ort1}.

\begin{lem}
\label{ort2}
Let $\delta_1,\delta_2\in \Delta_\el(G)$ and $\phi\in\Phi_2(G,\zeta)$. Then
\[
\int_{\Phi_2(G,\zeta)}n(\phi)S(\delta_2,\phi)\overline{S(\delta_1,\phi)}d\phi = n(\delta_1)\delta(\delta_1,\delta_2).
\]
\end{lem}

\begin{proof}
The proof is based on an application of the simple local stable trace formula from \cite[\S9-10]{ArtTW}. Let $f = f_1\times\bar{f}_2$ with $f_i\in \C_\cusp(G,\zeta)$, hence $f_{i,G}(\gamma)$ is supported on $\Gamma_\el(G)$ for $i=1,2$.  Since $G$ is quasisplit, we have a stable linear form $\C_\cusp(G,\zeta)$ given by
\[
S^G_\disc(f) = \int_{\Phi_2(G,\zeta)}n(\phi)^{-1}f^G_1(\phi)\overline{f^G_2(\phi)}d\phi
\]
that is equal to
\[
S^G(f)=\int_{\Delta_\el(G/Z)}n(\delta)^{-1}f^G_1(\delta)\overline{f^G_2(\delta)}d\delta.
\]
We first consider the spectral expansion. Applying the relation \eqref{adjoint} to \eqref{fgphi}, we may write
\[
f^G_i(\phi) = \int_{\Delta_\el(G/Z)}S(\phi,\delta)f^G_i(\delta)d\delta = n(\phi)\int_{\Delta_\el(G/Z)}n(\delta)^{-1} \overline{S(\delta,\phi)}f^G_i(\delta)d\delta.
\]
We then vary $f_i$ in a manner such that $f^G_i$ has compact support modulo $Z(F)$ on $\Delta(G)$ and so that $f_{i}^G$ approaches the $\zeta^{-1}$-equivariant Dirac measure at the image of $\delta_iZ(G)$ in $\Delta_\el(G)$ respectively for $i=1,2$. The function $f_i^G(\phi)$ thus approaches $n(\phi)n(\delta_i)^{-1}\overline{S(\delta_i,\phi)}$, and $S_\disc^G(f)$ approaches
\[
n(\delta_1)^{-1}n(\delta_2)^{-1}\int_{\Phi_2(G,\zeta)} n(\phi)S(\delta_2,\phi)\overline{S(\delta_1,\phi)} d\phi.
\]
On the geometric side, we see that as $f_i$ approach the Dirac measures on $\delta_iZ(F)$ respectively, the geometric expansion
\[
\int_{\Delta_\el(G/Z)}n(\delta)^{-1}f^G_1(\delta)\overline{f^G_2(\delta)}d\delta
\]
approaches $n(\delta_1)^{-1}\delta(\delta_1,\delta_2)$, and equating both sides, the identity follows. 
\end{proof}

We now state the orthogonality relation for the distributions $\Theta_{\tilde\xi'}(\delta',\delta)$ and $\Theta_{\tilde\xi'}^\mathrm{ad}(\delta',\delta)$. Note that our proof involves an interchange of the order of integration of the stable kernels, understood distributionally. We delay the proof of this property to \S\ref{proof}, where we shall require it again for \eqref{fubini}. 

\begin{prop}
\label{adj}
Given $\delta',\delta'_1\in \Delta^\F(G')$ for $G'\in \F_\el(G)$, we have
\be
\label{adj1}
\int_{\Delta(G/Z)}n(\delta)\Theta_{\tilde\xi'}^\mathrm{ad}(\delta',\delta)\Theta_{\tilde\xi'}(\delta'_1,\delta) d\delta =n(\delta')\tilde\delta(\delta',\delta'_1).
\ee
Similarly, given $\delta,\delta_1\in \Delta(G)$, we have
\be
\label{adj2}
\int_{\Delta(G')}n(\delta')\Theta_{\tilde\xi'}^\mathrm{ad}(\delta',\delta)\Theta_{\tilde\xi'}(\delta',\delta_1) d\delta' =n(\delta)\delta(\delta,\delta_1).
\ee
\end{prop}
\begin{proof}
The first identity will be a consequence of the second. We start with the first. We will first show that for any $\delta',\delta'_1\in \Delta^\F_\el(G')$, we have
\be
\label{dell1}
\int_{\Delta_\el(G/Z)}n(\delta)\Theta_{\tilde\xi'}^\mathrm{ad}(\delta',\delta)\Theta_{\tilde\xi'}(\delta'_1,\delta) d\delta =n(\delta')\tilde\delta(\delta',\delta'_1).
\ee
Then the required formula will follow from \eqref{extended} and the decomposition of the integral over $\Delta(G)$ into
\[
 \sum_{\{M\}}|W(M)|^{-1}\int_{\Delta_\el(M)}n(\delta_M)\Theta_{\tilde\xi'}^\mathrm{ad}(\delta',\delta_M)\Theta_{\tilde\xi'}(\delta'_1,\delta_M) d\delta_M,
\]
where we note that the sum contains at most one nonzero term.
From the definitions, we first write \eqref{dell1} as 
\[
\int_{\Delta_\el(G/Z)}n(\delta)\int_{\Phi(\tilde G',\tilde\zeta')}S'(\delta',\phi')S(\tilde\xi'\circ\phi',\delta)d\phi' \int_{\Phi(\tilde G',\tilde\zeta')}\overline{S'(\delta'_1,\phi'_1)S(\tilde\xi'\circ\phi'_1,\delta)}d\phi'_1  d\delta.
\]
Suppose first that the integration over $\Delta_\el(G/Z)$ and $\Phi(\tilde G',\tilde\zeta')$ can be interchanged. Then the integral over $\delta$ can be evaluated using the orthogonality relation \eqref{ort1} for $S(\phi,\delta)$ and \eqref{phim}. It follows then that the latter is equal to
\[
\int_{\Phi(\tilde G',\tilde\zeta')}\delta(\tilde\xi'\circ\phi',\tilde\xi'\circ\phi'_1)n(\tilde\xi'\circ\phi')S'(\delta',\phi')d\phi' \int_{\Phi(\tilde G',\tilde\zeta')}\overline{S'(\delta'_1,\phi'_1)}d\phi'_1,
\]
and reducing to the terms with $\phi=\phi_1$, the two integrals combine to 
\[
\int_{\Phi(\tilde G',\tilde\zeta')}n(\tilde\xi'\circ\phi')S'(\delta',\phi')\overline{S'(\delta'_1,\phi')}d\phi'.
\]
We claim that $n(\tilde\xi'\circ\phi')$ equals $n(\phi')$, so that the orthogonality relation from Lemma \ref{ort2} yields the required identity \eqref{dell1}. 

To prove the claim, let us compare 
\[
n(\tilde\xi'\circ\phi')= \int_{\Delta_\el(G/Z)}n(\delta)S(\tilde\xi'\circ\phi',\delta)\overline{S(\tilde\xi'\circ\phi',\delta)}d\delta
\]
with 
\[
n(\phi')= \int_{\Delta_\el(G')}n(\delta')S'(\phi',\delta')\overline{S'(\phi',\delta')}d\delta',
\]
where we note that the latter integrand depends only on the image of $\delta'\in \Delta_\el(\tilde{G}')$ in the set $\Delta_\el(\tilde{G}')/\tilde{Z}'(F)=\Delta_\el(\tilde{G}'/\tilde{Z}')=\Delta_\el(G')$. We can define an inner product on $\SI(G)$ by
\be
\label{innerp}
(a^G,b^G) = \int_{\Delta(G/Z)}n(\delta)^{-1}a^G(\delta)\overline{b^G(\delta)}d\delta,
\ee
whose restriction to $\SI_\cusp(G,\zeta)$ reduces to an integral over elliptic elements $\Delta_\el(G)$. Note that any function in $\SI(G)$ is bounded on $\Delta(G)$. Since the families of functions $\{n(\delta)\overline{S(\phi,\delta)}\}$ and $\{n(\delta')\overline{S'(\phi',\delta')}\}$ are orthogonal bases of $S\I_\cusp(G,\zeta)$ and $S\I_\cusp(\tilde{G}',\tilde\zeta')$ respectively, the identity will follow from showing that the stable transfer map is an isometry. Let us then consider 
\[
(a^{G'},b^{G'})= \int_{\Delta(G')}n(\delta')^{-1}a^{G'}(\delta')\overline{b^{G'}(\delta')}d\delta'.
\]
Once again the integrand is
\[
n(\delta')^{-1}\int_{\Delta(G/Z)}\int_{\Delta(G/Z)} \Theta_{\tilde\xi'}(\delta',\delta_1) \overline{\Theta_{\tilde\xi'}(\delta',\delta_2)}  a^G(\delta_1) \overline{b^G(\delta_2)}d\delta_1 d\delta_2,
\]
then by \eqref{adjdef} and \eqref{adj2} we see that the inner product is equal to
\[
 \int_{\Delta(G/Z)}\int_{\Delta(G/Z)} n(\delta_2)^{-1}\delta(\delta_1,\delta_2) a^G(\delta_1) \overline{b^G(\delta_2)}d\delta_1 d\delta_2,
\]
and evaluating at $\delta_1 = \delta_2$, we obtain $(a^G,b^G)$ as desired.

It remains to prove the second required identity \eqref{adj2}. In this case, beginning the argument parallel to the above leads to
\[
\int_{\Phi(\tilde G',\tilde\zeta')}\delta(\phi',\phi_1)n(\phi')S(\tilde\xi'\circ\phi',\delta)d\phi' \int_{\Phi(\tilde G',\tilde\zeta')}\overline{S(\tilde\xi'\circ\phi'_1,\delta_1)}d\phi'_1,
\]
and hence
\[
\int_{\Phi(\tilde G',\tilde\zeta')}n(\phi')S(\tilde\xi'\circ\phi',\delta)\overline{S(\tilde\xi'\circ\phi',\delta_1)}d\phi'.
\]
We see that this closely resembles the orthogonality relation of Lemma \ref{ort2}, and indeed we shall use a variation on the proof of the latter. In particular, recall that we may choose a suitable family of test functions $f_i\in\C_\cusp(G,\zeta)$ such that  $f_i^G(\phi)$ approaches $n(\phi)n(\delta_i)^{-1}\overline{S(\delta_i,\phi)}$ for $i=1,2$. Replacing $\phi$ by $\tilde\xi'\circ\phi'$ and choosing $\tilde\zeta'$ compatibly, we thus obtain a family of functions on $\C_\cusp(\tilde{G}',\tilde\zeta')$, which we write as $\tilde f_i$, so that the above equation is given as the limit of 
\[
\int_{\Phi(\tilde G',\tilde\zeta')}n(\phi')^{-1}f_1(\tilde\xi'\circ\phi')\overline{f_2(\tilde\xi'\circ\phi')}d\phi' = \int_{\Phi(\tilde G',\tilde\zeta')}n(\phi')^{-1}\tilde f_1(\phi')\overline{\tilde f_2(\phi')}d\phi' 
\]
as $f_1,f_2$ vary, where we note that the cuspidality of $\tilde f_i$ follows from that of $f_i$. That is, $f_i$ is supported on the set $\Gamma_{\el}(G)$ which we can identify with a subset of $\Gamma_{G,\el}(G')$ by Section \ref{ss}. Applying the local stable trace formula in this case to $\tilde{G}'$, we have that the latter is equal to
\[
\int_{\Delta_\el(G')}n(\delta')^{-1}\tilde f^G_1(\delta')\overline{\tilde f^G_2(\delta')}d\delta'.
\]
Since $f_i$ is chosen so that $f_{i}^G$ approaches the $\zeta^{-1}$-equivariant Dirac measure at the image of $\delta_iZ(G)$ in $\Delta_\el(G)$, it follows that $\tilde f_i$ vanishes unless $\delta'$ is an image of some $\delta_i$. Moreover, for such $\delta'$ we have 
\[
\pi_0((\hat{G}_{\delta'}')^\Gamma/Z(\hat{G}')^\Gamma)= \pi_0(\hat{G}_{\delta_i}^\Gamma/Z(\hat{G})^\Gamma),
\]
since $\hat{G}_{\delta'}'$ is isomorphic to $\hat{G}_{\delta_i}$ and $G'$ is elliptic, so that $n(\delta') = n(\delta_i)$,  thus giving $n(\delta)\delta(\delta_1,\delta_2)$. 
\end{proof}

\section{Stable transfer: Geometric and spectral}
\label{stransfer}

\subsection{The stable transfer conjecture}

Let us now return to our main Conjecture \ref{maincon}. We may reformulate it again as follows, with the proposed formula \eqref{candidate} for the transfer in Corollary \ref{skt}.

\begin{con}
\label{conjt}
For every $f \in \C(G,\zeta)$, there exists an $f'\in \SI(\tilde G',\tilde \zeta')$ such that
\begin{equation}
\label{strans}
f'(\delta') = \int_{\Delta(G/Z)} \Theta_{\tilde\xi'}(\delta',\delta)f^G(\delta) d\delta, \qquad \delta'\in\Delta_G(\tilde{G}')
\end{equation}
from $G$ to $G'$, where $f^G(\delta)$ denotes the stable orbital integral of $f$ at a strongly regular stable conjugacy class $\delta$. 
\end{con}

\noindent The stable transfer depends on the choice of auxiliary data, transfer factors, and Haar measures. In particular, we can view the conjecture as a transfer of Haar measures from $G$ to $G'$. We shall establish this at the end of the section.

As in the case of endoscopy, our transfer factors are defined only up to normalisation, so it is more appropriate to speak of families of transfer factors. We shall say that $f\in \C(G)$ and $f'\in\C(G')$ have matching (stable) orbital integrals if there exists a distribution $\Theta_{\tilde\xi'}(\delta',\delta)$ on $\Delta(G')\times \Delta(G)$ such that \eqref{strans} holds for all $\delta\in\Delta_G(\tilde{G}')$. Further, we call $\Theta_{\tilde\xi'}(\delta',\delta)$ a stable transfer factor if for each $f\in \C(G)$ there exists $f'\in\C(G)$ such that $f$ and $f'$ have matching orbital integrals.  We may as well require that $\Theta_{\tilde\xi'}(\delta',\delta)$ be nonzero only if $\delta'$ is an image of $\delta$. Conjecture \ref{conjt} then can be rephrased as the existence of a function $f' \in \C(\tilde G',\tilde\zeta')$ with matching stable orbital integrals and, implicitly, that our proposed distribution \eqref{candidate} is a stable transfer factor in the latter sense. 

Of course, we remind the reader that this conjecture is by no means new. We discuss some known or simple cases.

\ben
\item
When $G'=\{1\}$, it is trivially verified in \cite[p.178]{ST}. In that case, $f'$ is a constant, equal to the integral over $\Delta(G)$ of the product of $f^G(\delta)$ with the stable character $S(\phi,\delta)$, the latter being equal to $\Theta_{\tilde\xi'}(1,\delta)$. 

\item
When $G=\mathrm{SL}(2)$ and $G'$ a torus, this is again verified in \cite[\S2]{ST} and in the archimedean case, \cite[\S27]{shelstad} (see also \cite[\S8]{thomas}). We explain the transfer factors here in brief. In the split case $G' = \mathrm{GL}(1)$, the representations are simply one-dimensional characters $\chi$. The stable character on $G$ is a stably invariant function on regular semisimple elements, evaluating on the split torus to 
\[
|D^G(\delta(t))|^{-1}(\chi(t)+ \chi^{-1}(t)),
\]
where we embed $\GL(1)$ by the usual $\delta(t) =\text{diag}(t,t^{-1})$, and zero on elliptic classes. Then for $\delta\in G(F)$, the stable transfer factor can be computed by comparing stable characters on $G$ and $G'$,
\[
\Theta(\delta',\delta)  = |D^G(\delta)|^{-1}(\delta(\delta,\delta') + \delta({\delta^{-1}},\delta')),
\]
where $\delta(\cdot,\cdot)$ is understood as the delta distribution as in \cite[(2.13)]{ST}. 

In the nonsplit case, for the real elliptic torus $G'(\R) =T(\R)$, we parametrise its elements by $s(\theta)$ with $0\le \theta <2\pi$.  The stable transfer factor $\Theta_{\tilde\xi'}(s(\theta),\delta)$ is given by 
\[
\sum_{n\in\Z}e^{in\theta}\frac{\mp e^{in\theta}}{|e^{i\theta}-e^{-i\theta}|} \quad\text{or}\quad \sum_{n\in\Z}e^{in\theta}\frac{e^{nt} + e^{-nt}}{|e^t- e^{-t}|},
\]
depending on whether $\delta$ lies in the elliptic or the split torus of $G$ respectively. In the $p$-adic case, the cases separate into whether $G'$ is ramified, and in both cases the stable transfer is computed explicitly in \cite[\S2.4]{ST}, where we refer the reader for explicit formulas.

\een

\noindent Now let us examine slightly more general cases, {\em without} using the local Langlands correspondence. We call a function $ f\in \C(G,\zeta)$ cuspidal if $f_M$ vanishes for every proper Levi subgroup $M$ of $G$ (see also Section \ref{spaces}). We denote by $\C_\cusp(G,\zeta)$ the subspace of cuspidal functions. It is the subspace of $\C(G,\zeta)$ whose image in $\SI(G,\zeta)$ equals $\SI_\cusp(G,\zeta)$. 

Let $f\in \C_\cusp(G,\zeta)$  The property that $f$ is cuspidal implies that the image $f'$, if it exists, must vanish unless there are elliptic maximal tori $T\subset G$ and $T'\subset G'$ with admissible $L$-embeddings $^LT\subset {^LG}$ and $^LT'\subset {^LG}'$ such that $\xi'({^LT'})$ is contained in $^LT$. The problem thus reduces to that of tori. Similarly, it is also possible to consider minimal Levi subgroups $M\subset G$ and $M'\subset G'$, which are maximal tori, and restricting to stable conjugacy classes in $M(F)$ and $M'(F)$ respectively, though we will not study this here. 

\begin{lem}
\label{indaux}
The transfer $f\to f'$ is independent of choice of auxiliary datum.
\end{lem}

\begin{proof}
Suppose $(\tilde{G}'_1,\tilde{\xi}'_1)$ and $(\tilde{G}'_2,\tilde{\xi}'_2)$ are two auxiliary data with fixed central data $(Z_1,\zeta_1)$ and $(Z_2,\zeta_2)$, and associated stable transfer factors $\Theta_{\tilde\xi'_1}$ and $\Theta_{\tilde{\xi}'_2}$. Let $\tilde{G}'_{12}$ be the fibre product of $\tilde{G}'_1$ and $\tilde{G}'_2$ over $G'$. We have 
\[
Z(\hat{\tilde G}'_{12}) = (Z(\hat{\tilde{G}}'_1)\times Z(\hat{\tilde{G}}'_2))/\text{diag}_{-}(Z(\tilde{G}'))
\]
where $\text{diag}_{-}(Z(\tilde{G}'))$ is the anti-diagonal embedding. Given $w\in W_F$, let $g_w = (g(w),w)$ be an element in $\G'$ such that ad$_{g_w}$ acts by $\rho'$ on $\hat{G}'$. Also let 
\[
\tilde\xi'_i(g_w) = (\zeta_i(w),w),\qquad \zeta_i(w)\in Z(\hat{\tilde G}'_i)
\]
for $i=1,2$. Let $z_{12}(w)$ be the image of $(z_1(w),z_2(w)^{-1})$ in $Z(\hat{\tilde G}'_{12}) $, which is a cocycle of $W_F$ valued in $Z(\hat{\tilde G}'_{12})$, and by duality determines a character $\tilde\eta_{12}'$ of $\tilde G'_{12}$. Its restriction to $\tilde C_1\times \tilde C_2$ determining the fibre product $\tilde{G}'_1\times \tilde{G}'_2$ over ${G}'$ is equal to $\tilde\eta_1'\times (\tilde\eta'_2)^{-1}$, and pulling back the central datum $(Z,\zeta)$ we obtain the character $\tilde\zeta'_{12}$ on $\tilde{G}'_{12}$. For $i=1,2$, let $\delta_i'\in\Delta_G(\tilde{G}'_i)$ be 
 such that $(\delta_1',\delta_2') \in \tilde G'_{12}$. Then it follows that 
\be
\label{tf12}
\Theta_{\tilde\xi'_2}(\delta_2',\delta) = \tilde\zeta'_{12}(\delta_1',\delta_2') \Theta_{\tilde\xi'_1}(\delta_1',\delta).
\ee
from Lemma \ref{lem1}(i) and the definition of $\tilde\zeta'_{12}$.

The isomorphism $\SI(\tilde{G}'_1,\tilde\zeta_1)$ with $\SI(\tilde{G}'_2,\tilde\zeta_2)$ induced by the linear isomorphism $f_1\to f_2$ from $\C(\tilde{G}'_1,\tilde\zeta_1)$ to $\C(\tilde{G}'_2,\tilde\zeta_2)$ defined by 
\[
f_2(\delta_2') = \tilde\zeta_{12}(\delta_1',\delta_2')f_1(\delta_1'),
\]
where $\delta_1$ is any element such that $(\delta_1,\delta_2)\in \tilde{G}'_{12}.$ The isomorphism commutes with the transfer mappings $f\to f'_i = f^{\tilde{G}'_i}$. Then taking the inductive limit over such maps we see that the stable transfer mapping is independent of choice of auxiliary datum.
\end{proof}

\noindent As a consequence, the distribution $f'(\phi')$ depends on $\phi'$ rather than $\tilde\phi' = \tilde\xi'\circ\phi'$. However, $f'(\phi')$ does still depend on the choice of transfer factor. We shall think of $\Theta_{\tilde\xi'}$ as a family of transfer factors, one for each choice of $(\tilde{G}',\tilde\xi')$. Let us briefly indicate this. The (absolute) Langlands-Shelstad transfer factor that we have been discussing is based on the canonical relative transfer factor
\[
\Delta(\delta^e,\gamma,\bar\delta^e,\bar\gamma), \qquad \delta',\bar\delta'\in\Delta_{G}(\tilde G^e),\ \gamma,\bar\gamma\in\Gamma_{G}(G),
\]
associated to each $G,G^e,$ and $(\tilde G^e,\tilde\zeta^e)$ (see for example \cite[\S2]{ArtTW}). Recall that our assumption that $\tilde\xi'$ is of unitary type ensures that $|\Delta(\delta^e,\gamma,\bar\delta^e,\bar\gamma)|=1$. The pair $(\bar\delta^e,\bar\gamma)$ are chosen base points used to define the absolute transfer factor
\be
\label{abstf}
\Delta(\delta^e,\gamma) = \Delta(\delta^e,\gamma,\bar\delta^e,\bar\gamma)\Delta(\bar\delta^e,\bar\gamma),
\ee
defined to be zero unless $\delta^e$ is an image of $\gamma$. If we call an absolute transfer factor any function $\Delta(\delta^e,\gamma)$ on $\Delta_G(\tilde G')\times \Gamma(G)$ such that \eqref{abstf} holds if $\delta^e$ is an image of $\gamma$, and is zero otherwise, then the space of absolute transfer factors forms a $U(1)$-torsor. Following \cite[\S2]{Lnote}, we call $\Delta$ a transfer family for $(G,G^e)$ that varies according to $(\tilde G^e,\tilde\zeta^e)$, uniquely determined up to a multiplicative constant of absolute value one. 

As the stable transfer factor $\Theta_{\tilde\xi'}$ depends on the transfer family $\Delta$, we can in particular define a stable transfer family depending on $\Delta$. The relation \eqref{tf12} allows us to relate stable transfer factors associated to different auxiliary data. Similarly, suppose $t:G'_1\stackrel{\sim}{\to} G'_2$ is an isomorphism of mesoscopic data equipped with a dual $L$-isomorphism $\hat{t}:\G'_2\stackrel{\sim}{\to} \G'_1$. and $(\tilde{G}'_1,\tilde\xi'_1)$ is an auxiliary datum for $G'_1$, we obtain an auxiliary datum 
\[
(\tilde{G}'_2,\tilde\xi'_2) = t(\tilde{G}'_1,\tilde\xi'_1)
\]
 for $G'_2$ such that $t$ canonically extends to an $F$-isomorphism $\tilde{G}'_1\stackrel{\sim}{\to}\tilde{G}'_2$. We also have canonically a corresponding $L$-isomorphism $^Lt:{^L{\tilde{G}'_1}}\stackrel{\sim}{\to}^L{\tilde{G}'_2}$ such that 
 \[
 \tilde\xi'_2 = {^Lt}\circ\tilde\xi'_1\circ\hat{t}^{-1}.
 \]
 Then if $\Theta_{\tilde\xi'_1}(\delta'_1,\delta)$ is a transfer factor for $(\tilde{G}'_1,\tilde\xi'_1)$, it follows that  $\Theta_{\tilde\xi'_1}(t\delta'_1,\delta)$ is a transfer factor for $(\tilde{G}'_2,\tilde\xi'_2)$. A similar relation also holds in the endoscopic case.

As regards the general case, we have as usual the following reduction, which will allow us to work with $G'$ in place of $\tilde G'$ in many cases.

\begin{lem}
If Conjecture \ref{conjt} holds for $G$ with $G_\der$ simply connected and $(Z,\zeta)$ trivial, then it holds for arbitrary $G$ and $(Z,\zeta)$.
\end{lem}

\begin{proof}
Suppose that $G,G'$ and $(Z,\zeta)$ are arbitrary. Let $\tilde{G}$ be a $z$-extension of $G$ by the central induced torus $\tilde{C}$, and let $(\tilde Z,\tilde\zeta)$ be the pullback of $(Z,\zeta)$ to $\tilde{G}$. Since $G(F) = \tilde{G}(F)/\tilde{C}(F)$, we can identify $\I(\tilde{G},\tilde{\zeta})$ with $\I(G,\zeta)$. Moreover, there is a  bijection between isomorphism classes of mesoscopic data $\F(G)$ and $\F(\tilde{G})$. For any $G'$, we can find a extension $\tilde{\G}'$  equipped with an $L$-isomorphism $\tilde\xi': \tilde{\G}'\to {^L\tilde{G}'}$ and a natural embedding of $\G'$ into $\tilde{\G}'$, so that $(\tilde{G}',\tilde{\zeta}')$ is an auxiliary datum for both $G'$ and $\tilde{G}'$. We have a natural projection $\pi$ of $\I(\tilde{G})$ onto $\I(\tilde{G},\tilde{\zeta})=\I(G,\zeta)$ given by
\[
{a}({\gamma})\mapsto \int_{\tilde{Z}(F)}{a}(z{\gamma})\tilde{\zeta}(z)dz
\]
for any $a\in\I(\tilde{G})$. Similarly, we define a projection $\pi'$ of $\SI(\tilde{G}')$ onto $\SI(\tilde{G}',\tilde{\zeta}')=\SI({G}',{\zeta}')$. It then follows from Lemma \ref{lem1}(i) that the projections commute with the stable transfer mappings from $\I(\tilde{G})$ to $\SI(\tilde{G}')$ and $\I({G},\zeta)$ to $\SI({G}',\zeta')$ respectively. That is,
\begin{center}
\begin{tikzcd}
\I(\tilde{G}) \arrow[r] \arrow[d, "\pi"] & \SI(\tilde{G}') \arrow[d, "\pi'"] \\
\I(G,\zeta) \arrow[r]                  & \SI({G}',{\zeta}')
\end{tikzcd}
\end{center}
where the horizontal maps are the stable transfer mappings. The result follows.
\end{proof}

\subsection{Geometric transfer and spectral transfer}

The stable geometric transfer leads us naturally to a proof of the following stable character identity, which is simply a restatement of the conjectural formula \cite[(2.3)]{ST} in different terms.
\begin{cor}
\label{conj2}
The stable spectral transfer \eqref{bemap} implies that for any $\delta\in\Delta(G)$, we have 
\be
\label{schar}
S(\tilde\xi'\circ\phi',\delta) = \int_{\Delta(G')} \Theta_{\tilde\xi'}(\delta',\delta)S'(\phi',\delta')d\delta'.
\ee
\end{cor}

\begin{proof}
As noted in \cite[\S2.1]{ST}, this follows from the stable spectral transfer \eqref{bemap}, which follows by Propositions \ref{thm1} and \ref{skt}, and the fact that the functions $S(\phi,\delta)$ and $S'(\phi',\delta')$ are bases of the respective spaces $\SI_\cusp(G,\zeta)$ and $\SI_\cusp(\tilde G' ,\tilde \zeta')$ respectively \cite[Lemma 5.1]{LCR}, which then extends to the full spaces $\SI(G,\zeta)$ and $\SI(\tilde G' ,\tilde \zeta')$ as before.

We also note that this can also be seen directly from a combination of geometric and spectral transfer. Consider the stable spectral transfer (Proposition \ref{thm1}). Using the Fourier expansions for $f'(\phi')$ in \eqref{fgphi}, we have
\be
\label{specexp}
\int_{\Delta(G')}S'(\phi',\delta')f'(\delta')d\delta'  = \int_{\Delta(G/Z)}S(\tilde\xi'\circ\phi',\delta) f^G(\delta) d\delta.
\ee
Applying the stable geometric transfer, Corollary \ref{skt} which follows from \eqref{bemap}, the lefthand side is equal to
\be
\label{specexp2}
\int_{\Delta(G')}S'(\phi',\delta')\int_{\Delta(G/Z)} \Theta_{\tilde\xi'}(\delta',\delta)f^G(\delta) d\delta d\delta'.
\ee
Letting $f$ approach the $\zeta^{-1}$-equivariant Dirac measure at a fixed $\delta_1 Z(F)$, we have
\[
\int_{\Delta(G')}S'(\phi',\delta') \Theta_{\tilde\xi'}(\delta',\delta_1) d\delta' = S(\tilde\xi'\circ\phi',\delta_1)
\]
as required.
\end{proof}

In \cite{ST}, the integration is taken over the Steinberg-Hitchin base, the variety of stable semisimple conjugacy classes of $G'(F)$. Its measure is determined by the Haar measure on $G'$, in particular the singular locus has measure zero, and conincides with $\Delta(G')$.

\begin{rem} 
We can also give a heuristic derivation of the previous result using \eqref{candidate}. By definition, the right hand side of \eqref{schar} is 
\[
 \int_{\Delta(G')}  \int_{\Phi(\tilde G',\tilde\zeta')}S'(\delta',\phi'_1)S(\tilde\xi'\circ\phi'_1,\delta)d\phi'_1 S'(\phi',\delta')d\delta',
\]
and supposing we may interchange integrals, we then have
\[
 \int_{\Phi(\tilde G',\tilde\zeta')} \int_{\Delta(G')} S'(\delta',\phi'_1)S'(\phi',\delta')d\delta'  S(\tilde\xi'\circ\phi'_1,\delta)  d\phi'_1.
\]
Then applying Lemma \ref{adj0} to the inner integral leads to the orthogonality relation \eqref{ort1} for the stable virtual characters $S(\phi,\delta)$,
\[
\int_{\Delta_\el(G')} n(\delta')S(\phi',\delta')\overline{S'(\phi_1',\delta')}d\delta'  = \delta(\phi',\phi_1')n(\phi'),
\]
for a positive real number $n(\phi')$. Then evaluating the outer integral at $\phi'=\phi_1'$, we obtain $S(\tilde\xi'\circ\phi',\delta)$ as desired. 
\end{rem}

As was observed heuristically in \cite[\S26]{shelstad}, the stable transfer can be easily seen to be functorial in the following sense. Given $G$, let $G' \in \F(G)$ and $G''\in \F(G')$ with accompanying auxiliary data. 

\begin{cor} 
With hypotheses as above, we have $f^{\tilde G''} = (f^{\tilde G'})^{\tilde G''}$.
\end{cor}

\begin{proof}
The result can be seen to hold for stable characters by composing maps of Paley-Wiener spaces using the argument in Proposition \ref{thm1} and Corollary \ref{conj2}.
\end{proof}

\begin{rem}
To see this on the level of stable orbital integrals heursitically, assume for simplicity that $\G'$ and $\G''$ are $L$-groups, which we identify as $^LG'$ and $^LG''$ respectively. Then we have $L$-embeddings 
\[
^LG'' \stackrel{\xi''}{\to} {^LG'} \stackrel{\xi'}{\to} {^LG},
\]
and denote by $\xi'''$ the composition. First, we claim that
\be
\label{claim}
\Theta_{\xi'''}(\delta'',\delta)= \int_{\Delta(G')} \Theta_{\xi''}(\delta'',\delta')\Theta_{\tilde\xi'}(\delta',\delta)d\delta'.
\ee
To see this, we expand the transfer factor $\Theta_{\xi''}(\delta'',\delta')$ in the righthand side
\[
\int_{\Delta(G')}\int_{\Phi(\tilde G'',\tilde\zeta'')}S''(\delta'',\phi'')S'(\tilde\xi''\circ\phi'',\delta')d\phi''\Theta_{\tilde\xi'}(\delta',\delta)d\delta'.
\]
Formally interchanging the order of integration, this is
\[
\int_{\Phi(\tilde G'',\tilde\zeta'')}S''(\delta'',\phi'')\int_{\Delta(G')}S'(\tilde\xi''\circ\phi'',\delta')\Theta_{\tilde\xi'}(\delta',\delta)d\delta'd\phi'',
\] 
and we can apply Corollary \ref{conj2} to see that the inner integral is equal to $S(\tilde\xi' \circ \tilde\xi''\circ\phi'',\delta) = S(\tilde\xi'''\circ\phi'',\delta)$. In other words, we have
\[
\int_{\Phi(\tilde G'',\tilde\zeta'')}S''(\delta'',\phi'')S(\tilde\xi'''\circ\phi'',\delta)d\phi'',
\]
which is equal to $ \Theta_{\xi'''}(\delta'',\delta)$, the lefthand side.
Finally, writing $(f^{\tilde G'})^{\tilde G''}$ as 
\[
\int_{\Delta(G')} \Theta_{\xi''}(\delta'',\delta') \int_{\Delta(G/Z)} \Theta_{\tilde\xi'}(\delta',\delta)f^G(\delta) d\delta d\delta',
\]
then interchanging integrals and applying \eqref{claim} gives the result. 
\end{rem}

\subsection{Proof of Conjecture \ref{conjt}}
\label{proof}
By definition we may write \eqref{conjt} as 
\be
\label{exp}
\int_{\Delta(G/Z)}\int_{\Phi(\tilde G',\tilde\zeta')}S'(\delta',\phi')S(\tilde\xi'\circ\phi',\delta)d\phi'  f^G(\delta)d\delta.
\ee
On the other hand, Let $f'(\phi')$ be the transfer given by \eqref{bemap}. Applying the inversion formulas  together with spectral transfer, we have
\begin{align}
f'(\delta') &= \int_{\Phi(\tilde G',\tilde\zeta')}S'(\delta',\phi')f'(\phi')d\phi\notag\\
&= \int_{\Phi(\tilde G',\tilde\zeta')}S'(\delta',\phi')f(\tilde\xi'\circ\phi')d\phi\notag\\
&= \int_{\Phi(\tilde G',\tilde\zeta')}S'(\delta',\phi')\int_{\Delta(G/Z)}S(\tilde\xi'\circ\phi',\delta) f^G(\delta) d\delta\ d\phi, \label{fubini}
\end{align}
and in particular we note that this last expression is bounded. So we see that the result amounts to interchanging orders of integration. We shall consider this in two ways. 

Firstly, consider \eqref{fubini}. The stable kernels $S,S'$ are smooth and locally integrable on $G(F)$, and the inner integral converges absolutely for $f\in\C(G,\zeta)$. 
The outer integral, by \eqref{phi2} and \eqref{phim}, decomposes into
\[
\sum_{\{M'\}}|W(M')|^{-1}\sum_{\phi'\in \Phi_2(\tilde M',\tilde\zeta')/i\a^*_{\tilde{M}'}}\int_{i\a^*_{\tilde M',\phi'}}S'(\delta',\phi'_{\lambda'})S(\tilde\xi'\circ \phi'_{\lambda'},\delta)d\lambda',
\]
where in the integrand we have only indicated the terms depending on $\phi'$ for now.
Using the canonical projection of $i\a_M^*$ onto $i\a_G^*$, we can choose $\lambda\in i\a^*_{G}$ such that $(\tilde\xi'\circ\phi)_\lambda = \tilde\xi'\circ\phi_{\lambda'}$, since $i\a_{\tilde G'}^*$ acts on $\Phi(\tilde G',\tilde\zeta')$ and $i\a_{G}^*$ on $\Phi(G,\zeta)$ respectively, and the embedding $\tilde\xi'$ induces an embedding of $\Phi(\tilde G',\tilde\zeta')$ into $\Phi(G,\zeta)$. Moreover, recall from \cite[\S5]{LCR} the property that 
\[
S(\phi_\lambda,\delta) = S(\phi,\delta)e^{\lambda(H_G(\delta))}, \qquad \lambda\in i\a_{G}^*,\delta\in \Delta_\el(G).
\]
We see that by the adjoint relation \eqref{adjoint} and the identity $n(\phi_\lambda) = n(\phi)$ which follows from the definition of $n(\phi)$, the same property also holds for $S(\delta,\phi_\lambda)$. The latter integral is therefore equal to a sum over $M$ and $\phi'$ of 
\be
\label{conv}
S'(\delta',\phi')S(\tilde\xi'\circ \phi',\delta)\int_{i\a^*_{\tilde M',\phi'}}e^{\lambda'(H_{\tilde G'}(\delta'))}e^{\lambda(H_G(\delta))}d\lambda.
\ee
If $F$ is nonarchimedean, inner  integral is taken over compact tori $i\a^*_{\tilde M',\phi}$, while in the archimedean case, the absolute convergence follows from the growth conditions of the functions $I_M(\gamma,\pi)$ and $I_M(\pi,\gamma)$ from Theorems 4.1 and 4.3 of \cite{fourier}.  
In particular, the functions $S(\delta,\phi)$ and $S(\phi,\delta)$ are smooth and compactly supported on $\Phi(G,\zeta)$ when $F$ is nonarchimedean, and  
\[
|D_\lambda D_\delta S(\delta,\phi) | \le c(\delta) (1+ ||\mu_\phi||)^n
\]
where $n$ is a positive integer and $c(\delta)$ is a locally bounded function on $\Delta_\reg(G)$, both depending on a given pair of invariant differential operators $D_\delta$ and $D_\lambda$ transferred from $A_M(F)$ and $i\a^*_M$ respectively. Moreover, we write $\mu_\phi$ for the linear form that determines the infinitesimal character of $\phi$, which is a Weyl orbit of elements
in the dual of a complex Cartan subalgebra of $G$, equipped with a suitable Hermitian norm $||\cdot||$.   In particular, for every fixed $\phi'$, the given expression converges absolutely. 

So write \eqref{fubini} as 
\[
\lim_{t \to \infty }\sum_{\{M'\}}|W(M')|^{-1}\sum_{||\mu_{\phi'}|| \le t }\int_{i\a^*_{\tilde M',\phi'}}S'(\delta',\phi'_{\lambda'})\int_{\Delta(G/Z)}S(\tilde\xi'\circ \phi'_{\lambda'},\delta)f^G(\delta) d\delta\ d\lambda',
\]
where the inner sum runs over $\phi'\in \Phi_2(\tilde M',\tilde\zeta')/i\a^*_{\tilde{M}'}$ such that $||\mu_{\phi'}|| \le t$.   The integrands now converging absolutely, the Fubini theorem for measurable functions (e.g., \cite[Theorem 18.3]{fubini}) then allows for the change of order of integration. Finally, applying the dominated convergence theorem, we interchange the limit and integral to obtain the desired formula \eqref{conjt}.

Alternatively, define functionals on $\SI(G,\zeta)$ by
\[
T_{\phi'}:f^G(\delta)\to f^G(\tilde\xi'\circ\phi')  = \int_{\Delta(G/Z)}S(\tilde\xi'\circ\phi',\delta) f^G(\delta) d\delta 
\]
and on $\SI(\tilde G',\tilde \zeta')$ 
\[
T'_{\delta'}: f'(\phi') \to f'(\delta') =  \int_{\Phi(\tilde G',\tilde\zeta')}S'(\delta',\phi')f'(\phi')d\phi,
\]
which we may view as distributions on $\Delta(G/Z)$ and $\Phi(\tilde G',\tilde\zeta')$ respectively. Pre-composing with the mapping $f\to f'$ of \eqref{bemap} we can consider the latter also as a distribution on $G$. Then using the property that the tensor product  of distributions $T_{\phi'}\otimes T_{\delta'}'$ is commutative \cite[Theorem 40.4]{treves}, so the order of integration can be interchanged, interpreted distributionally.

\section{Stable transfer spaces and spectral transfer factors}
\label{spectf}

\subsection{Spaces of distributions}
\label{spaces}

Let $\D(G,\zeta)$ be the space of $\zeta$-equivariant invariant distributions that are supported on the preimage in $G(F)$ of finitely many conjugacy classes in $\bar{G}(F) = G(F)/Z(F)$. For $F$ nonarchimdean, it is equal to the space of ordinary orbital integrals, whereas if $F$ is archimedean, it also includes radial derivatives of orbital integrals.
 Let $\F(G,\zeta)$ be the space of $\zeta$-equivariant invariant distributions spanned by the invariant characters of $G(F)$, hence generated by characters attached to the set $\Pi(G,\zeta)$ of irreducible representations of $G(F)$ whose central character restricts to $\zeta$ on $Z(F)$. A distribution $D$ in either $\D(G,\zeta)$ and $\F(G,\zeta)$ can be regarded as a linear form 
\[
D(f) = f_G(D)
\]
on either $\H(G,\zeta)$ and $\I(G,\zeta)$.

Let $I$ be a continuous, invariant linear form on $\H(G,\zeta)$. We say that $I$ is supported on characters if $I(f)=0$ for any $f$ such that $f_G=0$. If so, then there is a continuous linear form $\hat{I}$ on $\I(G,\zeta)$ such that 
\[
\hat{I}(f_G)=I(f)
\]
for all $f\in \H(G,\zeta)$. If $D\in \F(G,\zeta)$, it is clear that it is supported on characters. On the other hand, if $D\in \D(G,\zeta)$, it can be expressed in terms of strongly regular invariant orbital integrals, which are supported on characters by \cite{itf2}. Together with the fact that characters are locally integrable functions, it follows that we can generate $\I(G,\zeta)$ by either irreducible tempered characters or strongly regular orbital integrals, both denoted $f_G$. 

We denote by $\I_\cusp(G,\zeta)$ the subspace of functions in $\I(G,\zeta)$ supported on $\Gamma_\el(G)$. If $Z$ contains the split component of the centre of $G$, there is a surjective linear map
\[
\F(G,\zeta)\to \I_\cusp(G,\zeta)
\]
canonically given by the elliptic virtual characters $I(\tau,\gamma)$ associated to any $D\in\F(G,\zeta)$. There is also a canonical linear section defined by the set $T_\el(G,\zeta)$ in $\F(G,\zeta)$ whose image forms a basis in $\I_\cusp(G,\zeta)$ \cite[\S4]{LCR}. We also denote by $\SI_\cusp(G,\zeta)$ the image of $\I_\cusp(G,\zeta)$ in $\SI(G,\zeta)$, and $\H_\cusp(G,\zeta)$ the preimage of $\I_\cusp(G,\zeta)$ in $\H(G,\zeta)$.

Let $S\D(G,\zeta)$ and $S\F(G,\zeta)$ be the stable subspaces of stable distributions in $\D(G,\zeta)$ and $\F(G,\zeta)$ respectively. Any distribution $S$ in $S\D(G,\zeta)$ and $S\F(G,\zeta)$ can be identified with a linear form $f^G\mapsto f^G(S)$ on $\SI(G,\zeta)$. We say a linear form $S$ on $\H(G,\zeta)$ is stable if its value at $f$ depends only on the endoscopic transfer $f^e$ in the case $G^e = G^*$. If $G$ is quasisplit, there is a unique linear form $\hat{S}$ on $\SI(G^*,\zeta^*)$ associated to $S$ such that 
\[
\hat{S}(f^*) = S(f)
\]
for any $f\in \H(G,\zeta)$. In general, for any stable distribution $S$ there is a unique continuous linear form $\hat{S}$ on $\SI(G,\zeta)$ such that
\[
\hat{S}(f^G) = S(f).
\]
One also has an alternative description of the cuspidal subspaces as follows. The restriction map $a^G\to a^M$ from $S\I(G(F))$ to $S\I(M(F))$ give a filtration 
\[
\F^M(S\I(G))=\{a^G\in S\I(G):a^L=0, L\subsetneq M\}
\]
of $S\I(G(F))$ over the partially ordered set $\L/W_0$. We can then identify $S\I_\cusp(G) = \F^G(S\I(G))$. Then the graded component
\be
\label{grading}
\G^M(S\I(G)) = \F^M(S\I(G))/\sum_{L\supsetneq M}\F^L(S\I(G))
\ee
attached to $\{M\}$ is canonically isomorphic to $S\I_\cusp(M)^{W(M)}$.

\subsection{Transfer spaces}

For any $G'\in \F_\el(G)$, let $\SI(\tilde{G}',G)$ be the subspace of functions in $\SI(\tilde{G}',\zeta)$ which depend only on the image of $\Delta_{G}(\tilde{G}')$ in $\tilde\Delta^\F(G)$, which we denote by $\Delta(\tilde{G}',G)$. We define the cuspidal subspace $\SI_\cusp(\tilde{G}',G) $ to be the intersection 
\[
\SI(\tilde{G}',G)\cap \SI_\cusp(\tilde{G}',\tilde\zeta')= \SI_\cusp(\tilde{G}',\tilde\zeta')^{\text{Out}_G(G')}.
\]
Assuming the stable transfer conjecture, it follows from the definitions that $f\to f'$ maps $\C(G)$ continuously to $\SI(\tilde{G}',G) $ and $\C_\cusp(G)$ continuously to $\SI_\cusp(\tilde{G}',G)$. If we define a function
\be
\label{a'}
a'(\delta') = a^{G'}(\delta') = \int_{\Delta(G/Z)} \Theta_{\tilde\xi'}(\delta',\delta) a^G(\delta)d\delta
\ee
on $\Delta_G(\tilde{G}')$, the stable transfer then gives a continuous map from $\SI(G)$ to $\SI(\tilde{G}',G)$ and $\SI_\cusp(G)$ to $\SI(\tilde{G}',G)$. Define the topological vector space
\[
S\I^\F_\cusp(G) = \bigoplus_{G'\in \F_\el(G)}\SI_\cusp(\tilde{G}',G)
\]
of smooth functions on $\Delta^\F_\el(G)$. For any function $a^G\in \SI_\cusp(G)$, we define the direct sum of images of $a^G$,
\[
a^\F=a^{G,\F} = \bigoplus_{G'\in \F_\el(G)}a'.
\]
Then the map 
\be
\label{TFmap}
\T^\F:a^G \to a^{G,\F}
\ee
is a continuous linear map from $\SI_\cusp(G)$ to $\SI^\F_\cusp(G)$. The following is a natural analogue of the endoscopic mapping in \cite[Proposition 3.5]{LCR} and \cite[I.4.11]{MW1} in the nonarchimedean case and \cite[I.4.12]{MW1} in the archimedean case. The most difficult part of lies in proving the surjectivity, which we shall return to in the next section. For now we simply take it on as an assumption. 

\begin{thm}
\label{TF}
Assume that $\T^\F$ is surjective. Then it is an isometric isomorphism.
\end{thm}

\begin{proof}
It is straightforward to see using the adjoint relations \eqref{adj1} and \eqref{adj2} that $\T^\F$ is invertible on its image, with inverse $a^{G,\F}\to a^G$ given by
\be
\label{inv}
a^G(\delta) = n(\delta)\int_{\Delta^\F_\el(G)}n(\delta')\Theta_{\tilde\xi'}(\delta,\delta')a^{G,\F}(\delta')d\delta'
\ee
for any $\delta\in\Delta(G)$ and $a^{G,\F}\in \SI^\F_\cusp(G)$. The map is moreover an isometry with respect to the inner product \eqref{innerp} on $\SI_\cusp(G)$ and
\be
\label{innerf}
(a^\F,b^\F) = \sum_{G'\in \F_\el(G)}\iota(G,G')(a',b')
\ee
on $\SI^\F_\cusp(G)$, where $\iota(G,G') = |\text{Out}_G(G')|^{-1}$. The inner product can first be expanded as 
\[
\sum_{G'\in \F_\el(G)}\iota(G,G')\int_{\Delta_{G,\el}(G')}n(\delta')^{-1} a'(\delta')\overline{b'(\delta')}d\delta',
\]
as $a',b'$ are cuspidal hence supported on the elliptic set. Expanding $a'$, the integrand is equal to
\[
n(\delta')^{-1} \int_{\Delta_\el(G/Z)} \Theta_{\tilde\xi'}(\delta',\delta) a^G(\delta)d\delta \overline{b'(\delta')}
= n(\delta') \int_{\Delta_\el(G/Z)}  a^G(\delta) \overline{ \Theta_{\tilde\xi'}(\delta,\delta') b'(\delta')} d\delta.
\]
Summing the integral over $G'$, we see that the constant $ |\text{Out}_G(G')|^{-1}$ normalises the measure on the quotient of $\Delta_{G,\el}(G')$ by $\text{Out}_G(G')$, then applying \eqref{inv} to $b^G$ we have
\begin{align*}
(a^\F,b^\F) &= \int_{\Delta^\F_\el(G)}  \int_{\Delta_\el(G/Z)} n(\delta')  a^G(\delta) \overline{ \Theta_{\tilde\xi'}(\delta,\delta') b^\F(\delta')} d\delta \ d\delta' \\
&= \int_{\Delta_\el(G/Z)} n(\delta)^{-1}a^G(\delta) \overline{b^G(\delta)}d\delta \\
&= (a^G,b^G)
\end{align*}
as required.
\end{proof}

\subsection{Surjectivity}

We now turn to the spectral analogue of our constructions so far. Our main goal will be to define stable spectral transfer factors. Assume that the spectral transfer \eqref{bemap} holds. We can then define a spectral basis parallel to $\Delta^\F_\el(G)$, namely, 
\[
\Phi^\F_2(G) = \coprod_{G'\in \F_\el(G)}\Phi_{2}(\tilde G',\tilde\zeta')/\text{Out}_G(G'),
\]
which can again be written as the set of pairs $(G',\phi')$. It parametrises a basis of $\SI^\F_\cusp(G)$. Also define
\be
\label{dec1}
\Phi^\F(G) = \coprod_{\{M\}}\Phi^\F_2(M)/W(M), 
\ee
which can also be described as the union over $W_0$-orbits $\{M\}$ in $\L$ and $W(M)$-orbits $\{M'\}$ in $\F_\el(M)$ of the quotient of $\Phi_2(\tilde M',\tilde\zeta')$ by Out$_M(M')\rtimes W(M)^{M'}$, where $W(M)^{M'}$ is the stabiliser of $M'$ in $W(M)$. We again have a decomposition according to central character,
\be
\label{dec2}
\Phi^\F(G) = \coprod_{\zeta}\Phi^\F(G,\zeta).
\ee
With these definitions in place, we return in earnest to the surjectivity of the map $\T^\F$ in \eqref{TFmap}, which is required in order to define our stable spectral transfer factors. It can be obtained as a consequence of the refined local Langlands correspondence.

\begin{lem}
\label{llcTF}
Suppose the refined local Langlands correspondence holds for $G$ over $F$.
Then $\T^\F$ is surjective.
\end{lem}
\begin{proof}
Again relying on the local Langlands correspondence over $F$, we can identify $\SI(\tilde G',\tilde \zeta')$ with the natural Schwartz space on a basis $\Phi(\tilde G',\tilde\zeta')$ of the vector space spanned by tempered, stable, $\tilde\zeta'$-equivariant characters on $\tilde G'(F)$, and similarly for $S\I(G,\zeta)$ and $\Phi(G,\zeta)$. The elements in $\Phi(\tilde G',\tilde\zeta')$ are indexed by tempered Langlands parameters $\phi'$ of $\tilde G'$, and in particular decomposes into a disjoint union of cuspidal Langlands parameters attached to Levi subgroups $\tilde M'$ of $\tilde G'$. We define a space $\tilde\Phi^\F(G,\zeta)$ analogous to $\tilde\Delta^\F(G)$ in Section \ref{conjcl}, which fibres over $\Phi^\F(G,\zeta)$, and $\SI^\F_\cusp(G,\zeta)$ can be identified with the natural equivariant Schwartz space on it as a consequence of the trace Paley-Wiener theorem for Schwartz functions on $G$ \cite{tracePW}.

Suppose that $\phi_1$ is a finite linear combination of linear forms in  $\tilde\Phi^\F(G,\zeta)$. We can assume that $\phi_1$  is the image of some $\phi'\in\Phi(\tilde G',\tilde\zeta')$ for some $G'\in \F_\el(G)$, such that $a^\F(\phi) = a'(\phi')$ for any $a^\F\in SI^\F(G,\zeta)$. The value at $\phi_1$ of any function $a^\F$ in $S\I^\F(G,\zeta)$ is then given by a finite linear combination
\[
a^\F(\phi_1) = \sum_{\phi'}c_{\phi'}a'(\phi'), \qquad \phi'\in\Phi(\tilde G',\tilde\zeta'), G'\in \F_\el(G).
\]
As an invariant distribution on $\tilde G'(F)$, any $\phi'$ can be identified with a locally integrable function whose restriction to $\Delta_G(\tilde G')$ is smooth. The set $\Delta_G(\tilde G')$ maps onto an open subset of $\tilde\Delta^\F(G)$ with finite fibres, and we can thus write $a^\F(\phi_1)$ as
\[
\int_{\Delta^\F(G/Z)} I^\F(\phi_1,\delta_1) a^\F(\delta_1)d\delta_1
\]
for some smooth, $\zeta$-equivariant function $I^\F(\phi_1)$ on $\tilde\Delta^\F(G)$, whose integral against any $a^\F=a^{G,\F}$ converges with respect to the measure $d\delta_1$. Applying \eqref{a'}, we have then
\[
\int_{\Delta(G/Z)}  I(\phi_1,\delta) a^G(\delta)d\delta, \qquad a^G\in S\I(G,\zeta),
\]
where
\[
I(\phi_1,\delta) = \int_{\Delta^\F(G/Z)} I^\F(\phi_1,\delta_1) \Theta_{\tilde\xi'}(\delta_1,\delta)d\delta_1
\]
is again a smooth, $\zeta$-equivariant function on $\Delta(G)$, whose integral against any $a^G$ converges with respect to the measure $d\delta$. Letting $a^G$ now approximate the $\zeta^{-1}$-equivariant Dirac measure at $\delta Z(F)$, it follows that if the function
\[
a \to a^\F(\phi_1), \qquad a\in \C(G,\zeta)
\]
induced by $\phi_1$ vanishes, then so does $I(\phi_1,\delta)$. By the inversion formula \eqref{adj1}, it follows that $I^\F(\phi_1,\delta_1)$ also vanishes on $\tilde\Delta^\F(G)$, and hence $\phi_1$ itself vanishes. The mapping $\T^\F$ is thus locally surjective.

To see that it is surjective, we can define a spectral transfer factor describing the local mapping $a\to a^\F$,
\[
a^\F(\phi_1) = a'(\phi') = \int_{\Phi_2(G,\zeta)} \Theta_{\tilde\xi'}(\phi_1,\phi)a(\phi)d\phi, \qquad a\in\C(G,\zeta),
\]
compatible with the decompositions \eqref{dec1} and \eqref{dec2} as below, hence with the characterisations of $S\I(G,\zeta)$ and $S\I^\F(G,\zeta)$ as Schwartz spaces of functions on $\Phi(G,\zeta)$ and $\tilde\Phi^\F(G,\zeta)$ respectively. Thus the surjectivity extends.
\end{proof}

\subsection{Stable spectral transfer factors}

Our discussion here now parallels the geometric case, so we may be brief. Let $f \in \C_\cusp(G,\zeta)$. For any $G'\in\F_\el(G)$, the transfer $f'$ is a function in $\SI_\cusp(\tilde{G}',\tilde\zeta')^{\text{Out}_G(G')}$ and $f'(\phi')$ is defined for every $\phi'\in\Phi_2(\tilde{G}',\tilde\zeta')$. 
If $\T^\F$ is surjective, we may define a stable spectral transfer factor $\Theta_{\tilde\xi'}(\phi',\phi)$ to be any distribution on $\Phi_2(\tilde{G}',\tilde\zeta')\times \Phi_2(G,\zeta)$ such that the identity
\[
f'(\phi') = \int_{\Phi_2(G,\zeta)}\Theta_{\tilde\xi'}(\phi',\phi)f(\phi)d\phi
\]
holds. The stable transfer factors can again be extended to distributions on $\Phi^\F(G)\times \Phi(G)$ as follows. Set $\Theta_{\tilde\xi'}(\delta',\delta)$ to be zero unless there is an $M$ such that $(\phi',\phi)$ belongs to the Cartesian product of $\tilde\Delta^\F_{G,\el}(M)/W(M)$ with $\Delta_{G,\el}(M)/W(M)$. If there is such an $M$, then $(\phi',\phi)$ is the image of a pair $(\phi'_M,\phi_M)$ in $\Phi^\F_{2}(M)\times \Phi_{2}(M)$, and we set
\be
\label{extended2}
\Theta_{\tilde\xi'}(\phi',\phi) = \Theta_{\tilde\xi',G}(\phi',\phi) =\sum_{w\in W(M)}\Theta_{\tilde\xi',M}(\phi'_M,w\phi_M),
\ee
where again each sum contains at most one nonzero term, and depends only on $\phi'$ and $\phi$.  Finally, define the adjoint spectral transfer factor 
\[
\Theta_{\tilde\xi'}(\phi,\phi') =n(\phi')^{-2}\overline{\Theta_{\tilde\xi'}(\phi',\phi)},
\]
which complements the adjoint stable geometric transfer factors quite nicely. As in the geometric case, their definition is imposed upon us by the adjoint relations they satisfy, parallel to Proposition \ref{adj}.

\begin{prop}
Given $\phi',\phi'_1\in \Phi^\F(G')$ for $G'\in \F_\el(G)$, we have
\be
\label{adj3}
\int_{\Phi(G)}n(\phi)\Theta_{\tilde\xi'}(\phi',\phi)\Theta_{\tilde\xi'}(\phi,\phi'_1) d\phi =n(\phi')\tilde\delta(\phi',\phi'_1).
\ee
Similarly, given $\phi,\phi_1\in \Phi(G)$, we have
\be
\label{adj4}
\int_{\Phi(G')}n(\phi')\Theta_{\tilde\xi'}(\phi',\phi)\Theta_{\tilde\xi'}(\phi_1,\phi') d\phi' =n(\phi)\delta(\phi,\phi_1).
\ee
\end{prop}
\begin{proof}
Since we do not have an explicit description of the spectral transfer factors, the proof will follow instead from interpreting the linear isometry $\T^\F$ spectrally. We recall the spectral form of the inner product on $\SI_\cusp(G)$,
\[
(a^G,b^G) = \int_{\Phi_2(G)}n(\phi)^{-1}a^G(\phi) \overline{b^G(\phi)}d\phi, \qquad a^G,b^G\in\SI_\cusp(G).
\]
We then have the spectral form of the inner product on $\SI^\F_\cusp(G)$ in \eqref{innerf},
\[
(a^\F,b^\F) = \sum_{G'\in \F_\el(G)}\iota(G,G')\int_{\Phi_2(G')}n(\phi')^{-1}a'(\phi') \overline{b'(\phi')}d\phi'.
\]
Then defining the spectral analogue of the inverse of $\T^\F$,
\[
a^G(\phi) = n(\phi)\int_{\Phi^\F_2(G)}n(\phi')\Theta_{\tilde\xi'}(\phi,\phi')a^{G,\F}(\phi')d\phi',
\]
and using the definition of $\Theta_{\tilde\xi'}(\phi',\phi)$ and its adjoint above, it follows by the same argument as in the proof of Proposition \ref{TF} that $(a^\F,b^\F) = (a^G,b^G)$. In particular, $\Theta_{\tilde\xi'}(\phi',\phi)$ and $\Theta_{\tilde\xi'}(\phi,\phi')$ represent kernels of inverse transforms of each other. 
\end{proof}


\appendix

\section{On the surjectivity of $\T^\F$ (without local Langlands)}
\label{appb}

The proofs of surjectivity of the analogous endoscopic map for nonarchimedean fields \cite[Lemma 3.4]{LCR} and \cite[I.4.11]{MW1} both involve a reduction to the Lie algebra at the identity element, either by passing to germs of orbital integrals or Harish-Chandra descent to the Lie algebra of unipotent subgroups. In both cases, the descent of transfer factors and properties of the Fourier transform on the Lie algebra are employed. In order to discuss the latter, we need some preparations.

\subsection{Descent principles}
\label{descsec}
We explore the relationship between descent and transfer. The descent properties will be required later to define our stable spectral transfer factors. As in the case of endoscopy \cite{LS2}, this reduces to the descent of transfer factors, which we shall have to take on as a hypothesis in this paper. 
Given $G'\in\F(G)$ and $d'\in\Delta_\ss(G')$, let  $\tilde{d}'$ be its preimage in $\Delta_\ss(\tilde{G}')$. We shall investigate the behaviour of 
\[
f'(\tilde\delta') = \int_{\Delta(G/Z)} \Theta_{\tilde\xi'}(\tilde\delta',\delta)f^G(\delta) d\delta
\]
for $\tilde\delta'$ near to $\tilde d'$. If $d'$ is not the image of any semisimple element in $\Delta_\ss(G)$, then no strongly $G$-regular element in $G'_{d'}(F)$ can be the image of an element in $G(F)$, thus $f'$ vanishes on $\Delta_G(\tilde{G}'_{\tilde{d}'})$. It follows then that $f'$ vanishes for all $\tilde{d}'$ in a neighbourhood of $\tilde{d}'$ in $\tilde{G}'(F)$. We shall therefore assume that $d'$ is the image of some $d\in \Delta(G)$. 

We can find a representative $d'_1 = x^{-1}d' x$ in the stable conjugacy class of $d'$ such that $G'_{d'_1}$ is quasisplit over $F$, and we can multiply $x$ by an element of $G'_{d'_1}$ if necessary so that 
\[
\text{Int}(x^{-1}): G'_{d'} \to G'_{d'_1}
\]
is defined over $F$. Then $x$ acts on the preimage $\tilde G'_{d'}$ of $G'_{d'}$ in $\tilde G'$, so that 
$
f'(x^{-1}\tilde\delta' x) = f'(\tilde\delta') 
$
for all $\tilde\delta'\in\Delta_G(\tilde{G}_{d'})$, which can also be seen to follow from the same property for stable orbital integrals \cite[\S1.3]{LS2}. In particular, replacing $d'$ by $d'_1$ if necessary, we may assume that $G'_{d'}$ is quasisplit over $F$.

The group $G'_{d'}$ induces a mesoscopic datum for $G_d$ in the following sense. More generally, we call $d'$ a $T'$-image of $d$ in $G$ if there exists an admissible embedding of tori $T'\to T^*$ sending $d'$ to $d^*$ in $G^*$ and an $x\in G^*_\text{sc}$ such that 
\[
(\text{Int}(x)\circ \psi) (d) = d^*
\]
and both $\text{Int}(x)\circ \psi$ and the preimage of $T$ are defined over $F$. Varying over $T'$ we obtain all images of $d$. Let $d'$ be a $T'$-image of $d$ for some torus $T'$, and let $d^*$ be the image of $d'$ under an admissible embedding of $T'$ in $T^*$, which we may choose to be such that $G^*_{d^*}$ is quasisplit. Fixing $G'\in \F(G)$, we shall attach an extension $\G_{d'}'$ of $W_F$ by $\hat{G}'_{d'}$ and an admissible embedding $\xi_{d'}': \G'_{d'}\to {^LG'_d}$ such that $(G'_{d'},\G_{d'}',\xi_{d'}')$ is a mesoscopic datum for $G_d$ and $T'\to T^*$ is admissible.

\begin{lem}
\label{embed}
The embedding $T'\to T^*$ can be chosen to be admissible for both $(G,G')$ and $(G_d,G'_{d'})$, unique up to isomorphism. Any admissible embedding of a maximal torus of $G'_{d'}$ in $G^*_{d^*}$ is admissible as an embedding of a maximal torus of $G'$ in $G^*$ and sends $d'\to d^*$. 
\end{lem}

\begin{proof}
The proof is a simple modification of \cite[1.4]{LS2}. We supply the details here for the sake of completeness. We first explain how the mesoscopic data is constructed. Given an embedding $T'\to T$, let $B'$ and $B^*$ be the associated Borel subgroups and $x\in G^*_\text{sc}$ as above. The map $\psi_x = \text{Int}(x)\circ\psi$ defines the quasisplit inner twist $G^*_{d^*}$ of $G_{d}$. The embedding $T' \to T^* \stackrel{\psi_x}{\leftarrow }T$ is dual to the diagram
\be
\label{diag}
\hat{G}'\leftarrow T'_1 \simeq T_1 \to \hat{T}^* \to \hat{T},
\ee
by which we can identify the set of coroots $R(G,T)^\vee$ of $T$ in $G$ with the set of roots $R(\hat{G},T_1)$ of $T_1$ in $\hat{G}$, and hence $R(G_{d},T)^\vee$ with a subset of $R(\hat{G},T_1)$.

Fix an $L$-group data, meaning a complex reductive group $\hat{G}_d$, an action $\rho_d$ of $\Gamma$ on $\hat{G}_d$, and a $\Gamma$-stable bijection $\Psi(G_d)^\vee\to \Psi(\hat{G}_d)$, by which we define $^LG_{d} = \hat{G}_d \rtimes W_F$. We may assume that $\hat{G}_d$ contains $T_1$ and that $R(\hat{G}_d,T_1)$ is equal to $R(G_d,T)^\vee$ as subsets of $R(\hat{G},T_1)$. Let $B_d^* = B^*\cap G^*_{d^*}$, and let $\mathcal B_d$ be the Borel subgroup of $\hat{G}_d$ generated by $T_1$ and the $B_1$-positive roots of $T_1$ in $\hat{G}_d$. We can then identify the map $T_1 \to \hat{T}^*$ in \eqref{diag} with the embedding $\hat{T}^*\to T_1$ in $\hat{G}_d$ given by $B_d^*$ and $\mathcal B_d$. The isomorphism 
\[
\hat{T}\stackrel{\psi_x}{\to} \hat{T}^*\to T_1
\]
yields an embedding of $\hat{T}$ in $\hat{G}_d$ and extends to an admissible embedding of $^LT$ in $^LG_d$, whose image is independent of the choice of extension.

Given $G'$, the dual $\hat{G}'_{d'}$ of $G'_{d'}$ is a subgroup of $\hat{G}_d$ normalised by $^LT$. We define $\G'_{d'}$ to be the subgroup of $^LG$ generated by $\hat{G}'_{d'}$ and $^LT$, and note that it is contained in $\G'$. Then we can also define $\xi'_{d'}$ to be the map given by restriction of $\xi'$ from $\G'$ to $\G'_{d'}$. We have a split exact sequence 
\[
1 \to \hat{G}'_{d'}\to \G'_{d'}\to W_F\to 1,
\]
and it follows that $(G'_{d'},\G_{d'}',\xi_{d'}')$ is a mesoscopic datum for $G_d$. Finally, we identify the embedding $\hat{T}'\to T_1$ given by $B'\cap G'_{d'}$ and $\mathcal B_{d}\cap \hat{G}'_{d'}$ with the restriction of the embedding $\hat{T}\to T_1$ above. This gives the admissible embedding $T'\to T^*$ as desired.

The choice of $\hat{G}_{d}$ is unique up to isomorphism of mesoscopic data. First suppose $B,B'$ are changed but $T'\to T^*$ remains fixed. Then the $L$-data $(\hat{G}_d,\rho_d)$ is replaced by another pair $(\hat{G}_d^1,\rho_d^1)$ that is $\Gamma$-isomorphic to it that sends the root datum $R(\hat{G}_d,T_1)$ to $R(\hat{G}^1_d,T_1)$, the image of $^LT$ in $^LG_d$ to its image in $^LG^1_d$, and $\hat{G}'_{d'}$ and $\G_d'$ to the new $\hat{G}^{1'}_{d'}$ and $\G_d^{1'}$ respectively. In particular, this gives an isomorphic mesoscopic datum for $G_d$.

If we replace $T'\to T^*$ with another admissible embedding of tori $T^{1}{}'\to T^{1*}$, such that $d'$ lies in $T',T^1{}'$ and $d^*$ in $T^*,T^{1*}$. Then we may assume that the new Borel subgroups are obtained from $B'$ and $B^*$ by conjugation in $G'_{d'}$ and $G^*_{d^*}$ respectively, and the new data is isomorphic. Finally, it is straightforward to see that the choice of $G' \in \F(G)$ within its isomorphism class does not affect the isomorphism class of $G'_{d'} \in \F(G_d).$ 
\end{proof}


We shall say that $(G,G')$ admits $\Theta_{\tilde\xi'}$-transfer if for each $f\in C_c^\infty(G(F))$ there exists $f'  = f^{\tilde{G}'} \in C_c^\infty(\tilde{G}'(F))$ such that $f,f'$ have $\Theta_{\tilde\xi'}$-matching orbital integrals, that is,
\[
f'(\tilde\delta') = \int_{\Delta(G/Z)} \Theta_{\tilde\xi'}(\tilde\delta',\delta)f^G(\delta) d\delta.
\]
As we are not in a position to prove the descent of stable transfer factors, we simply admit it as a hypothesis.

\begin{lem}
\label{desclem}
Suppose that there exists some constant $c$ such that
\be
\label{tfdesc}
\Theta_{\tilde\xi'}(\tilde\delta',\delta)  \Theta_{\tilde\xi'_d}(\tilde\delta',\delta)^{-1} \to c
\ee
as $\tilde\delta\to \tilde d'$ and $\delta\to d$. Then $f'(\tilde\delta')$ is equal to a finite linear combination of stable orbital integrals on $G'_{d'}$. 
\end{lem}
\begin{proof}
According to the measure on $\Delta(G)$, the integral decomposes into
\[
\sum_{\{M\}}|W(M)|^{-1}\sum_{\{T\}}|W_F(G,T)|^{-1}\int_{T(F)}\Theta_{\tilde\xi'}(\tilde\delta',t)f^G(t) dt,
\]
where the inner sum is over stable conjugacy classes of elliptic maximal tori of $G$ over $F$, and $W_F(G,T)$ is the subgroup of elements in the absolute Weyl group of $(G,T)$ defined over $F$. Moreover, our hypothesis \eqref{tfdesc} implies 
\[
\int_{T(F)}\Theta_{\tilde\xi'}(\tilde\delta',t)f^G(t) dt = c\int_{T(F)}\Theta_{\xi'_d}(\tilde\delta',t)f^G(t) dt.
\]
Now by Lemma \ref{embed}, we may choose an embedding $T'\to T^*$ that is admissible for both $(G,G')$ and $(G_d,G'_{d'})$. If $\tilde\delta'$ is an element in the preimage of $T'(F)$ in $\Delta_G(\tilde{G}')$, then the integral in $f'(\tilde\delta')$ is taken over the stable conjugacy classes defined by the composition
\[
\tilde{T}'\to T' \to T^* \to T
\]
in $\Delta(G)$. Note that it is possible that $\tilde\delta'$ is a $T'$-image for more than one tori, and varying over equivalence classes tori in $\tilde G'_d$ we obtain all possible images. Then using the property that $\Theta_{\xi'_d}(\tilde\delta',t)$ vanishes unless $\tilde\delta'$ is an image, we conclude that the righthand side can be written as a finite linear combination $f^{\tilde G'_{d'}}(\tilde\delta')$.
\end{proof}

It follows from Lemma \ref{lem1}(i) that for $\tilde\delta'$ close to the identity, $\Theta_{\tilde\xi'}(\tilde\delta',\delta)$ depends only on the image $\delta$ of $\delta'$ in $\Delta(G')$, so we write it as $\Theta_{\tilde\xi'}^\text{loc}(\delta',\delta)$. 
We say that $(G,G')$ admits local $\Theta_{\tilde\xi'}$-transfer at the identity if for any $f\in C_c^\infty(G(F))$  we have
\be
\label{loctrans}
f'(\delta')  = \int_{\Delta(G/Z)}\Theta_{\tilde\xi'}^\text{loc}(\delta',\delta)f^G(\delta) d\delta
\ee
for all $\delta' \in \Delta_G(\tilde{G}')$ near to the identity.

\begin{cor}
Let $F$ be nonarchimedean, and assume \eqref{tfdesc} holds. If $(G_d,G'_{d'})$ have local $\Theta_{\xi'_{d'}}$-transfer at the identity for all $d\in\Delta(G)$, then $(G,G')$ has $\Theta_{\tilde\xi'}$-transfer.
\end{cor}

\begin{proof}
Since the assumption continues to hold if $G$ is replaced by a $z$-extension $\tilde{G}$, we can assume that $G = \tilde G$ and $\G'$ is an $L$-group. By \cite[Lemma 2.2A]{LS2} it suffices to show that $f'(\delta')$ is a local stable orbital integral on $G'(F)$, in the sense that for every semisimple element $d$ in $G(F)$ there exists $f_d\in C_c^\infty(G(F))$ such that $f'(\delta) = f'_d(\delta)$ for all regular semimsimple $\delta$ near to $d$. Again by Lemma \ref{lem1}(i), it follows that local $\Theta_{\xi'_{d'}}$-transfer at the identity implies local $\Theta_{\xi'_d}$-transfer at $d'$ in the sense that \eqref{loctrans} holds for all $\delta' \in \Delta_G(\tilde{G}')$ near to $d'$. By assumption, we may apply the descent formula to express $f'(\tilde\delta')$ as a finite linear combination of stable orbital integrals on $G'_{d'}$. Then applying local transfer at $d'$ to each summand, the result follows.
\end{proof}


\subsection{A stable kernel identity}

Let $\g$ be the Lie algebra of $G$, and similarly $\g'$ of $G'$. Fix a symmetric, nondegenerate $G$-invariant bilinear form $B$ on $\g$ and a nontrivial additive character $\psi_0$ on $F$. For any $\varphi\in C_c^\infty(\g(F))$, we define the Fourier transform
\[
\hat \varphi(Y) = \int_{\g(F)}\varphi(X)\psi_0(B(X,Y))dX,
\]
which acts as a linear isomorphism from $C_c^\infty(\g(F))$ to itself. It a well-known result of Harish-Chandra that there exists a smooth, locally-integrable function
\[
i: \Gamma(\g)\times \Gamma(\g) \to \CC,
\]
where $\Gamma(\g)=\Gamma_\reg(\g(F))$ is the space of regular $G(F)$-orbits in $\g(F)$, such that 
\[
\varphi_G(X) = \int_{\Gamma(\g)}i(X,Y) (\hat \varphi)_G(Y)dY, \qquad X\in \Gamma(\g)
\]
for a fixed Haar measure $dY$ on $\Gamma(\g)$. As with $\Gamma(G)$, we may decompose the integral into a sum of integrals over conjugacy classes of maximal tori in $G$. Taking $G$ to be quasisplit, we define the smooth function
\[
s(S,T) = |\K_T|^{-1}\sum_{X\to S}\sum_{Y\to T}i(X,Y),
\]
where $S,T$ are regular stable $G(F)$-orbits in $\g_\reg(F)$, the sums run over the distinct $G(F)$-orbits in each respective stable orbit, and $|\K_T|$ is equal to the number of $Y$ in the orbit of $T$. We also have the stable analogue
\[
\varphi^G(S) = \int_{\Delta(\g)}s(S,T) (\hat \varphi)^G(T)dT, \qquad S\in \Delta(\g)
\]
which is a consequence of \cite[(3.3)]{LCR}. We define the Lie algebra analogue of the stable transfer factor 
\[
\Theta_{\xi'}(S',T) =\int_{\Delta_G(\g')}s'(S',T')s(d\xi'(T'),T)dT', \quad S'\in \Delta_G(\g'), T\in \Delta(\g),
\]
where $s'$ denotes the function associated to $G'$, and $d\xi'$ is the induced map on regular semisimple elements from $\g(F)$ to $\g'(F)$. Then suppose the following stable analogue of Waldspurger's kernel formula \cite[1.2]{transfert} holds for any $G'\in \F_\el(G)$,
\be 
\label{wkernel}
\int_{\Delta(\g)} \Theta_{\xi'}(S',S)s(S,T) dS = \delta_0 \int_{\Delta_G(\g')} s'(S',T')\Theta_{\xi'}(T',T)dT',
\ee
for $T\in \Delta(\g')$ and $S'\in \Delta_G(\g')$, where the latter is the set of regular stable $G(F)$-orbits in $\g'_\reg(F)$ and $\delta_0$ is a constant depending only on $G,G'$.  For any $\varphi\in C_c^\infty(\g(F))$, we define the transfer
\[
\varphi'(S') = \int_{\Delta(\g)}\Theta_{\xi'}(S',S)\varphi^G(S)dS,
\]
and by the preceding formulas, it is equal to
\begin{align}
\varphi'(S') &=\int_{\Delta(\g)}\Theta_{\xi'}(S',S) \int_{\Delta(\g)}s(S,T) (\hat{\varphi})^G(T)dT\notag\\
&= \delta_0 \int_{\Delta(\g)}\int_{\Delta_G(\g')} s'(S',T')\Theta_{\xi'}(T',T)(\hat{\varphi})^G(T)dT'dT\notag\\
&=  \delta_0 \int_{\Delta_G(\g')} s'(S',T') (\hat \varphi)'(T')dT'. \label{s'}
\end{align}

\subsection{Surjectivity}

With these considerations, the surjectvitity of the map $\T^\F$ in the nonarchimedean case, without local Langlands, can then be shown to follow from the proposed kernel formula and the descent of transfer factors.

\begin{lem}
Let $F$ be a non-archimedean local field, and assume \eqref{tfdesc} and \eqref{wkernel}. Then $\T^\F$ is surjective.
\end{lem}
\begin{proof}
Assume first that restriction of $\T^\F$ to $S\I_\cusp(G,\zeta)$ maps onto the corresponding cuspidal subspace $S\I_\cusp^\F(G,\zeta)$ of $S\I^\F(G,\zeta)$. Recall that the filtration on $\SI(G,\zeta)$ with respect to $\L/W_0$ gives a grading \eqref{grading} of the space, whereby
\[
S\I(G,\zeta) = \bigoplus_{\{M\}}S\I_\cusp(M,\zeta)^{W(M)},
\]
and similarly
\[
S\I^\F(G,\zeta) = \bigoplus_{\{M\}}S\I_\cusp^\F(M,\zeta)^{W(M)}.
\]
The map $\T^\F$ is compatible with these gradings, and the transfer mapping can identified with the corresponding transfer mapping for cuspidal functions on each $M$.

The surjectivity of $\T^\F$ on $S\I_\cusp(G,\zeta)$ will follow from the corresponding surjectivity of germs on the Lie algebra. Assume for simplicity that $G'$ is an $L$-group, hence $\tilde G' = G'$, and moreover that the central datum $(Z,\zeta)$ is trivial. Given $a^\F \in S\I_\cusp^\F(G)$, by the relation \eqref{adj2} we see that the function
\[
a^G(\delta) = \int_{\Delta^\F_\el(G)}\Theta_{\xi'}(\delta,\delta')a^\F(\delta')d\delta', \qquad \delta\in\Delta_\el(G)
\]
implies 
\[
a^\F(\delta') = \int_{\Delta_\el(G)}\Theta_{\xi'}(\delta',\delta)a^G(\delta)d\delta = (\T^\F(a^G))(\delta').
\]
Together with \eqref{adj1}, we see that $a^\F$ lies in the image of $S\I_\cusp(G)$ if and only if $a^G \in S\I_\cusp(G)$. We may assume that the components $a'$ of $a^\F$ are nonzero for exactly one $G'\in\F_\el(G)$. It suffices then to show that for $a^\F \in S\I_\cusp^\F(G)$, the function 
\be
\label{agint}
a^G(\delta) = \int_{\Delta_\el(G')}\Theta_{\xi'}(\delta,\delta')a'(\delta')d\delta', \qquad \delta\in\Delta_\el(G)
\ee
lies in $S\I_\cusp(G)$. Let $d$ be a fixed elliptic semisimple conjugacy class in $G(F)$, and let $S\G_\cusp(G,d)$ be the space of germs of functions in $S\I_\cusp(G)$ around $d$. Let $d'_1,\dots,d'_n$ be representatives of Out$_G(G')$-orbits of stable conjugacy classes in $G'(F)$ that are images of $d$, chosen such that $G'_{d'_j}$ is quasisplit. The image of $a^G$ in $S\G_\cusp(G,d)$ depends only on the image of $a'$ in the spaces of germs of functions in $S\I_\cusp(G')$ around $d'_j$. We may assume that these images are nonzero for exactly one $j$, so that the integral \eqref{agint} is supported on points close to $d_j'$. Then for any $\delta,\delta'$ close to $d,d'$ respectively and assuming the descent of transfer factors, Lemma \ref{desclem} reduces $a^G(\delta)$ and $a'(\delta')$ to orbital integrals on $G_d$ and. $G'_{d'}$. In particular, $d$ is central in $G_d$, and by Lemma \ref{lem1}(i) it suffices to show that $a^G(\delta)$ lies in $S\G_\cusp(G,1)$ for $\delta$ close to 1. 
 
Let $C_{c,\cusp}^\infty(\g(F))$ be the subspace of cuspidal functions on $\g(F)$, and let $S\G_\cusp(\g)$ be the space of germs of stable orbital integrals of such functions around 0. It is a finite dimensional space of germs of functions on $\Delta(\g)$. The exponential map gives a linear bijection from $S\G_\cusp(\g)$ to $S\G_\cusp(G,1)$. (This passage from orbital integrals on the group to germs on the Lie algebra can also be alternately described by Harish-Chandra descent, c.f. \cite[I.4.1]{MW1}.) Define also the finite dimensional vector space
\[
S\G_\cusp^\F(\g) = \bigoplus_{G'\in \F_\el(G)}S\G_\cusp(\g')^{\text{Out}_G(G')},
\]
and let 
\[
\tau^\F: S\G_\cusp(\g) \to S\G_\cusp^\F(\g) 
\]
be the Lie algebra analogue of $\T^\F$ on germs. We would like to show that it is surjective. Fix $G'\in\F_{\el}(G)$ and consider the image of an arbitrary $g'(S')\in S\G_\cusp(\g')^{\text{Out}_G(G')}$ in $S\G_\cusp^\F(\g)$, where $S'$ is belongs to the $G$-regular elliptic quotient $\Delta_{G,\el}(\g')/\text{Out}_G(G')$. In particular, $g'$ is supported on the regular elliptic locus. By \eqref{s'}, we can express $g'(S')$ as a finite linear combination
\[
g'(S') = \sum_{i=1}^n c_i s'(S',T'_i),\qquad T_i\in\Delta_{G,\el}(\g').
\]
We may choose the bilinear form $B'$ on $\g'$ to be invariant under Out$_G(G')$, so that the coefficients $c_i$ are constant on Out$_G(G')$ orbits. By Howe's finiteness theorem applied to $\g'(F)$ \cite[Theorem 2]{howe}, we can choose a compact neighbourhood $V_i$ of $T_i$ in $\Delta_{G,\el}(\g')$ for each $i$, such that $s'(S',T') = s'(S',T'_i)$ for all $T'\in V_i$ and $S'$ sufficiently small. Altogether, this implies that we can choose a $C_c^\infty$-function $\alpha'$ on $\Delta_{G,\el}(\g')^{\text{Out}_G(G')}$ such that
\[
g'(S') = \delta_0\int_{\Delta_{G,\el(\g')}}s'(S',T')\alpha'(T')dT'
\]
for all $S'$ sufficiently close to 0. The adjoint relations of Proposition \ref{adj} have Lie algebra analogues. Recall that the proof in the group case relies on the stable local trace formula, and similarly one may employ a stable local trace formula for the Lie algebra, which is simpler, to deduce the necessary relations. This allows us to invert the map 
\[
C_c^\infty(\Delta_\el(\g)) \to \bigoplus_{G_1'\in \F_\el(G)}C_c^\infty(\Delta_{G,\el}(\g')),
\]
thereby giving a function $\varphi_0\in C_c^\infty(\g_{\reg,\el}(F))$ such that for any $G_1'\in \F_\el(G)$, we have $\varphi_0^{G_1'}$ equals $\alpha'$ if $G_1'=G'$ and is trivial otherwise. Extending by zero, we can find a function $\varphi \in C_c^\infty(\g(F))$ such that $\hat\varphi = \varphi_0$. By \eqref{s'}, it follows that 
\[
\varphi'(S') = \varphi^{G'}(S') = g'(S'),
\]
and $\varphi^{G'_1}=0$ for any $G_1'\neq G'$. Thus the map $\tau^\F$ is surjective.
\end{proof}

\subsection*{Acknowledgments} The author was partially supported by NSF grant DMS-2212924. I am grateful to the anonymous reviewer for generous feedback, also to Julia Gordon and Daniel Johnstone for helpful discussions.

\bibliography{BESL2}
\bibliographystyle{alpha}
\end{document}